%% file: main.tex
\documentclass[11pt]{article}

\usepackage[a4paper,margin=1in]{geometry}

\usepackage[utf8]{inputenc}

\usepackage[ngerman,british]{babel}

\usepackage[T1]{fontenc}
\usepackage{mathpazo,tgcursor,tgpagella}
\usepackage{textcomp}
\usepackage{lmodern} 
\usepackage[english=british]{csquotes}

\usepackage{authblk}

\usepackage[bookmarksnumbered,bookmarksopen,unicode]{hyperref}
\hypersetup{
	pdftitle={},
	pdfauthor={},
	pdfkeywords={},
	pdfsubject={} 
}

\usepackage[textsize=tiny]{todonotes}
\usepackage{amsmath}
\usepackage{amssymb}
\usepackage{mathtools}
\usepackage[inline]{enumitem}
\usepackage{dsfont}
\usepackage{caption}
\usepackage{subfigure}
\usepackage{hyperref}
\usepackage[capitalize,noabbrev]{cleveref}
\usepackage{algorithm}
\usepackage{algorithmic}
\usepackage{booktabs}

\newcommand{\N}{\mathds{N}}
\newcommand{\Z}{\mathds{Z}}

\newcommand{\R}{\mathds{R}}

\newcommand{\OO}{\mathcal{O}}

\newcommand{\ra}{\rightarrow}

\newcommand{\abs}[1]{\lvert #1 \rvert}
\newcommand{\norm}[1]{\lVert #1 \rVert}

\newcommand{\iset}[1]{\mathcal{#1}}

\newcommand{\iS}{\iset{S}}

\newcommand{\iX}{\iset{X}}

\newcommand{\iP}{\iset{P}}

\newcommand{\set}[1]{\{#1\}}
\newcommand{\lrset}[1]{\left\{#1\right\}}

\newcommand{\qm}[1]{``#1''}

\newcommand{\Todo}[1]{\ifthenelse{\boolean{ToDoSwitch}}{\textcolor{red}{\textbf{ToDo:} #1 \\}}{}}

\newcommand{\ie}{i.e.\ }

\newcommand{\eg}{e.g.\ }

\newcommand{\st}{\mathrm{s.t.}}

\newtheorem{theorem}{Theorem}[section]
\newtheorem{corollary}[theorem]{Corollary}

\newtheorem{lemma}[theorem]{Lemma}

\newcommand{\scenario}{s}
\newcommand{\scenarioset}{\iS}

\newcommand{\mS}{\mathcal{S}}
\newcommand{\mP}{\mathcal{P}}
\newcommand{\mX}{\mathcal{X}}
\newcommand{\mI}{\mathcal{I}}
\newcommand{\Prob}{\mathbb{P}}

\usepackage{natbib}
\bibpunct[, ]{(}{)}{,}{a}{}{,}%
\title{Data-driven Distributionally Robust Optimization over Time}

\author[1]{Kevin-Martin Aigner\footnote{Corresponding author. E-mail address: \href{mailto:kevin-martin.aigner@fau.de}{kevin-martin.aigner@fau.de}}}
\author[1]{Andreas Bärmann} 
\author[1]{Kristin Braun}
\author[1]{Frauke Liers}
\author[2]{Sebastian Pokutta}
\author[3]{Oskar Schneider}
\author[2]{Kartikey Sharma}
\author[1]{Sebastian Tschuppik}
\affil[1]{Friedrich-Alexander University Erlangen-Nürnberg, Erlangen, Germany}
\affil[2]{Zuse Institute Berlin, Berlin, Germany}
\affil[3]{Fraunhofer Institute for Integrated Circuits IIS, Nürnberg, Germany}
\vspace{\baselineskip}

\newcommand{\proof}[1]{\emph{#1}}
\newcommand{\Halmos}{\hfill$\square$}

\begin{document}
\maketitle


\input{abstract}


\input{introduction}

\input{related_work}
\input{data_driven_dro_v2}

\input{dro_via_online_learning}
\input{numerical_results}

\input{conclusion}

\section{Acknowledgments}
We are grateful to Daniela Bernhard for proofreading parts of the manuscript.
We also thank the DFG for their support within Projects B06 and
B10 in CRC TRR 154, as well as within Project-ID 416229255 - SFB 1411.
This work has been supported by grant 03EI1036A from the Federal Ministry for Economic Affairs and Energy, Germany.


\bibliographystyle{apalike}



\bibliography{References} 



\clearpage
\input{extension}

\end{document}

%% file: abstract.tex
Stochastic Optimization (SO) is a classical approach
for optimization under uncertainty
that typically requires knowledge
about the probability distribution of uncertain parameters. 
As the latter is often unknown,
Distributionally Robust Optimization (DRO) provides a strong alternative
that determines the best guaranteed solution
over a set of distributions (ambiguity set).
In this work, we present an approach for DRO over time that uses online learning and scenario observations arriving as a data stream
to learn more about the uncertainty. 
Our robust solutions adapt over time
and reduce the cost of protection with shrinking ambiguity. 
For various kinds of ambiguity sets,
the robust solutions converge to the SO solution. 
Our algorithm achieves the optimization and learning goals
without solving the DRO problem exactly at any step.
We also provide a regret bound for the quality of the online strategy
which converges at a rate of $ \OO(\log T / \sqrt{T})$,
where $T$ is the number of iterations.
Furthermore, we illustrate the effectiveness of our procedure
by numerical experiments on mixed-integer optimization instances
from popular benchmark libraries
and give practical examples
stemming from telecommunications and routing.
Our algorithm is able to solve
the DRO over time problem significantly faster
than standard reformulations.

%% file: introduction.tex
\section{Introduction}

Many practical optimization problems
deal with uncertainties in the input parameters,
and it is important to compute optima that are protected against them.
Two prime methodologies to deal with uncertainty in optimization problems
are \emph{Stochastic Optimization (SO)} and \emph{Robust Optimization (RO)}. 
SO considers all uncertain parameters to be random variables,
and its solution approaches usually rely on the knowledge
of the probability distribution.
Classically, SO aims to find solutions that are optimal in expectation
(or more generally with respect to chance constraints
or different risk-measures, see~\cite{BirgeLouveaux1997}).
RO is typically used when knowledge about the probability distribution
is not at hand or a better guarantee of feasibility
is desired \citep{Ben-TalElGhaouiNemirovski2009}.
It strives to find solutions which perform best
against adversarial realizations of the uncertain parameters from a predefined uncertainty set.

Even if the underlying probability distributions are not at hand,
they can often be estimated from historical data.
These estimates are naturally also affected by uncertainty.
As such, recent research has focused on compromising between SO and RO
in order to obtain better protection against uncertainty
while controlling the conservatism of robust solutions.
In particular, \emph{Distributionally Robust Optimization (DRO)} aims to solve
a \qm{robust} version of a stochastic optimization problem
under incomplete knowledge about the probability distribution.
The benefit of DRO is that the solutions are fully protected
against the uncertainty and thus outperform non-robust solutions
with respect to worst-case performance.

Current research in DRO is primarily aimed at developing efficient solution techniques for static DRO problems, 
where the optimization problem is solved for a given and fixed ambiguity set.
However, in many practical applications,
additional information about the uncertainty becomes available over time.
For example, in situations where planning processes need to be repeated over time,
each plan may want to incorporate the outcomes of the previous decisions.
Applications abound for such processes, for example in
taxi or in ambulance planning. Applications also occur in iterative assigning 
landing time windows to aircraft such that security distances are kept
at an airport at all times even in case of disturbances which may lead
to frequent reassignments. For more details on the air traffic
application, see \citep{ATM1, Kapolke2016}.
Naturally, it would be beneficial to include any new data into a DRO approach
as soon as it becomes available.

In this article, we present a DRO approach
that iteratively incorporates such information over time.
Specifically, we provide an online learning algorithm
that solves DRO problems with limited initial knowledge about the uncertainty,
but which can leverage additional incoming data.
This allows the optimal solutions to adapt to the uncertainty
and gradually reduce the cost of protection.
To this end, we use scenario observations arriving as a data stream
to construct and update the ambiguity sets.
These sets contain the true data generating distribution
with high confidence and converge to it over time.
We also show that the solution to the DRO problem
converges to the true SO solution,
as the ambiguity sets shrink to the true distribution,
and hence the online algorithm also converges to the SO solution. 
However, the primary goal is to use the online algorithm
to solve the DRO problem.

The main feature of our work is an integrated procedure that can iteratively solve the DRO problem while learning reliable and time-dependent ambiguity sets.
We show that our online algorithm outperforms prior methods.
We also compare different approaches
to construct data-driven ambiguity sets.
In computational experiments,
we evaluate our algorithm on state-of-the-art benchmark libraries
and realistic case studies.
We demonstrate that our online method
leads to significantly reduced computation times
with only marginal sacrifices in solution quality.

\subsubsection*{Problem Statement.}

We consider the problem of minimizing a function 
$ f\colon \mathcal{X}\times \mathcal{S} \rightarrow \R $
over some (possibly non-linear and/or mixed-integer) feasible set~$ \mathcal{X} $.
We focus on the case of finitely many scenarios which are contained
in the set $ \mathcal{S} \coloneqq \set{s_1, \ldots, s_{\abs{\mathcal{S}}}} $.
Note that this is a natural modeling assumption that appears in several applications as many real-world random events
are best represented via discrete scenarios.
However, our approach is also able to treat the case of continuous random variables by sampling a sufficient large discrete scenario set from the probability distribution.
Finite approximations such as the sample average approximation~\citep{sample-average-approximation} and other similar scenario reduction techniques are standard in stochastic optimization.
They lead to algorithmic tractability in more general settings,
which is necessary for any realistically sized problem.
Furthermore, the use of a finite number of scenarios
allows us to treat the probability distribution as a vector
and hence leverage methods from first-order optimization
to solve the problem of finding the optimal distribution.
The setting of infinite-dimensional ambiguity sets for continuous distributions goes beyond the focus of this paper.

With this in mind, we define $ \mathcal{P}_0 \coloneqq
	\set{p \in [0, 1]^{\abs{\mathcal{S}}} \mid
		\sum_{k = 1}^{\abs{\mathcal{S}}} p_k = 1} $
as the $ \abs{\mathcal{S}} $-dimensional probability simplex.
We start with the probability simplex as our initial ambiguity set
as it imposes no restrictions on the distributions. 
However, we are not restricted to this choice and one can initialize with sets constructed using already available historical data. 
This would lead to less conservative solutions in practice without changing our theoretical results fundamentally.
Each point $ p \in \mathcal{P}_0 $ represents a probability distribution
over the scenarios $ s \in \mathcal{S} $.
Given a probability distribution of scenarios $ p^* \in \iP_0 $,
one can solve the following SO problem:
\begin{equation}
	\label{Eq:Stoch_nominal}
	\tag{SO}
	J^* \coloneqq \min_{x \in \iX} \mathbb{E}_{s \sim p^*}[f(x, s)]
		= \textstyle\sum_{k = 1}^{\mathcal{\iS}} f(x, s_k) p_k^*.
\end{equation}
However, if there is limited information about the probability vector~$ p^* $,
we can limit the impact of uncertainty
by solving the distributionally robust counterpart of the SO problem, namely
\begin{align}
	\min_{x \in \mathcal{X}} \; \max_{p \in \mathcal{P}_0}
	\; \mathbb{E}_{s \thicksim p}\left[f(x, s)\right]
	= \min_{x \in \mathcal{X}} \; \max_{s \in \mathcal{S}} \; f(x, s). \tag{DRO}\label{Eq:DRO_full}
\end{align}
This equality holds because the worst-case probability
is realized by some unit-vector
\linebreak
$ (0, \ldots, 0, 1, 0, \ldots, 0) \in \mathcal{P}_0 $.
For this maximally large ambiguity set,
the solutions to the DRO problem may be overly conservative.  
Furthermore, there is typically data available
in the form of observed realizations of the uncertain parameters.
The primary goal of this work is to develop an online algorithm
to solve the above DRO problem over time
while progressively integrating additional information.
This is achieved by refining the ambiguity sets over time
as we learn more about the uncertainty, e.g., with new realizations.
The use of DRO along with online optimization limits the impact
of adverse realizations.
Simultaneously, learning ensures that we can adapt and increase
our confidence as we gather more information.

We assume that information about the probability distribution $p^*$
is revealed in the form of i.i.d.\ realizations over time.
As such, we solve a sequence of DRO problems
with progressively shrinking ambiguity sets $\iP_t$.
At each time step $ t = 1, \ldots, T $ over a given horizon,
we construct the ambiguity set $\iP_t$
according to the confidence regions
estimated from scenario observations up to time~$t$.
Using these sets, we solve
\begin{align}
	\label{Eq:DRON}
	\tag{DRO$_t$}
	\widehat{J}_t := \min_{x \in \mathcal{X}} \;\; \max_{p \in \mathcal{P}_t}
	\;\; \mathbb{E}_{s \thicksim p}\left[f(x, s)\right].
\end{align}
In order to solve problem~\ref{Eq:DRON} efficiently,
our online optimization approach alternates
between a gradient step (for $p$)
and solving the minimization problem (for~$x$).
The task is to find a solution~$ x_t $ for each round $t$.
This approach can also be interpreted as a game over $T$ rounds,
where a player tries to make optimal decisions
against an adversary who chooses the probability distribution
from which to draw the uncertain realization.

In our analysis, we calculate the gap
between the worst-case performance of~$ x_t $ (generated by the online algorithm) over the set $ \iP_t $
and the performance of the optimal DRO solution at time~$t$. 
We prove that the average gap (reminiscent of the notion of dynamic regret)
is bounded as
\begin{alignat*}2
	&\frac{1}{T} \sum_{t=1}^{T} \Big(  \underbrace{\max_{p \in \mP_{t}}\mathbb{E}_{s \thicksim p}\left[f(x_{t},s)\right]}_{\text{ our online cost at time $t$}}  - \underbrace{\min_{x\in \mathcal{X}} \max_{p\in \mathcal{P}_{t}} \mathbb{E}_{s \thicksim p}\left[f(x,s)\right]}_{\text{optimal cost in hindsight at $t$}} \Big) 
	\leq \OO\Big(\sqrt{\frac{ h(T)}{T}} \Big),
\end{alignat*}
with high probability.
Here, $ h(T) $ is a bound on the path length of the distributions, i.e., $ \sum_{t = 1}^{T} \frac12 \norm{p_t - q_t}_2^2 \leq h(T) $
for $ p_t \in \mP_{t - 1} $ and $ q_t \in \mP_t $. 
We show that $ h(T) \in \OO(\log^2 T) $
for the different categories of ambiguity sets that we consider.
This bound controls the difference between the performance
of our online method and an exact DRO solver. 
It certifies that our approach successfully approximates
the DRO solution with an average gap that decreases over time.
We also show that the DRO solution converges to the true SO solution as the ambiguity sets converge to the true distribution in the limit.

%% file: related_work.tex
\subsubsection*{Related Work.}

Current research in distributional robustness is mostly concerned
with the appropriate choice of ambiguity sets to obtain guarantees on solution quality 
\citep{DelageYe2010,ParysEtAl2017}.
It also focuses on the derivation of algorithmically tractable reformulations
for the resulting robust counterparts, see \eg \cite{WiesemannKuhnSim2014},
as well as \cite{Calafiore2006}.
Ambiguity sets can be constructed by imposing constraints on expectation~\citep{ChenEtAL2017},
covariance~\citep{DelageYe2010}, mode~\citep{hanasusanto2015distributionally}, etc.
Another option is to use distance metrics
such as the Wasserstein metric~\citep{esfahani2018data},
$\phi$-divergence~\citep{bayraksan2015data},
$f$-divergence \citep{duchi2021statistics},
kernel-based distances~\citep{KirschnerBogunovicJegelkaKrause2020},
hypothesis tests~\citep{BertsimasEtAl2018a}, etc. 



Ambiguity sets can also be defined
by confidence bounds~\citep{RahimianMehrotra19},
and many popular ambiguity sets~\citep{BertsimasEtAl2018a, KirschnerBogunovicJegelkaKrause2020}
have associated probabilities of containing the true distribution. 
In this work, we construct the sets
as the combination of a simplex and a bounded set
defined by either confidence intervals, $\ell_2$-norm or kernel based metrics.
These bounded sets function as confidence regions.  

Online learning is an established field
which provides algorithms for solving problems over time. 
For a broad introduction see~\citet{hazan_ioco}.
Recently, this approach has been leveraged
to solve robust optimization problems~\citep{ben-tal:oracle_ro,ho2018online,pokutta2021adversaries,oracle_robust_optimization_nonconvex_loss}. 
Online learning has also been applied to DRO problems
where the ambiguity set is constructed from data. 
\citet{namkoong2016stochastic} leverage online optimization
for a DRO problem with ambiguity sets
defined by $f$-divergences.  
They use an alternating mirror descent algorithm
and provide regret bounds for the same. 
The authors of~\cite{online_method_dro_with_kl_regularization} propose a duality-free online stochastic method for a class of DRO problems with KL-divergence regularization on the dual variables.
Different algorithms are proposed and analyzed in~\cite{online_learning_for_convex_dro} for DRO problems with convex objectives with conditional value at risk and $\chi^2$-divergence ambiguity sets.
These existing works combine DRO and online learning in order to solve a single DRO problem.
This paper looks at a different problem.
The key difference between these works and ours
is that we consider a planning problem that has to be solved repeatedly in time with growing knowledge about the uncertain parameters.
Therefore, we establish an online framework to solve a series of DRO problem with changing ambiguity sets as new information arrives. 
This allows us to obtain robust online solutions
while learning about the uncertainty.


In \citet{KirschnerBogunovicJegelkaKrause2020},
the authors also focus on a DRO problem in an online learning context.
They allow for uncertainty in the parameters and only have noisy blackbox access to the objective function.
To solve the DRO problem while learning the objective and the ambiguity set, the authors solve a large program with convex constraints at each stage.
This was extended to more general ambiguity sets constructed with $\phi$-divergences in~\cite{drbo_2022}. 
In contrast, our work primarily focuses on obtaining distributionally robust optimal solutions in an online fashion while learning the true distribution. 
As a key advantage we point out that our algorithm does not require the
solution of the entire DRO problem at any step but only computes
gradient steps in the ambiguity space. 
Furthermore, we consider multiple ambiguity sets and leverage the online gradient descent algorithm to allow for faster computation and better applicability in real-world settings.

Another similar work
is by \citet{SessaBogunovicKamgarpourKrause2020}.
Therein, the authors present an online learning approach with a multiplicative weight algorithm in order to compute strict robust mixed-strategies over a decision set. 
In contrast to our work however, they do not consider a DRO setting. 
Finally, \citet{BaermannMartinPokuttaSchneider2020} and \citet{BaermannPokuttaSchneider2017}
consider the problem of learning unknown objectives in an online fashion but without any uncertainty in the models.


\subsubsection*{Contribution.}
The two key differences in our work which distinguish it from regular online optimization are
(i) use of DRO while learning from data
and (ii) solving the DRO problem approximately.
The first ensures that our solutions are robust to uncertainty
in the knowledge of the true probability distribution.
The second allows us to obtain robust solutions
without solving the DRO problem exactly at each step.
Specifically, the key contributions of our work are:
\\
\emph{Online Learning Algorithm for DRO}. We provide an online algorithm to solve the DRO  problem. 
It also learns the uncertainty from scenario observations over time,
shrinking the ambiguity sets. 
This allows for rapid computation of the DRO solutions along with their adaptation.
Thus, reducing the cost of protection. \\
\emph{Stochastic Consistency}. We also prove that the solution of the DRO problem
 converges to the SO problem.  
Since our online algorithm solves the DRO problem, thus it too converges to the solution of the SO problem. \\ 
\emph{High Probability Regret Bounds}. We prove that the cumulative regret between the solutions generated by our online method and the exact DRO solution at each time step shrinks at a rate of  $\mathcal{O}(\log T / \sqrt{T})$ with high probability.\\
\emph{Flexibility of Uncertainty Models}. We consider 3 different ambiguity models:  (i) confidence intervals, (ii) $\ell_2$-norm sets and (iii) kernel based ambiguity sets. 
These allow our approach to adapt to the application.\\
	%
\emph{Computational Results}. We provide a computational study on
mixed-integer benchmark instances and on real world problem
examples. Specifically, we compare on the MIPLIB and QPLIB libraries and further illustrate our results with two realistic applications from telecommunications and routing. 
In both cases our approach is considerably faster than solving the full distributionally robust counterparts. 

\subsubsection*{Outline}
Section~\ref{sec:dddro} presents the foundations of \emph{Data-Driven DRO}.
We present our algorithms and theoretical results on \emph{DRO over Time} in Section~\ref{sec:drotime}. 
Finally, in Section~\ref{sec:num_res},
we evaluate our methods on benchmark instances.

%% file: data_driven_dro_v2.tex
\section{Data-Driven DRO}
\label{sec:dddro}

In this section, we introduce ambiguity sets which form a key part of DRO problems.
We also introduce the dual reformulation
which is a standard way for solving robust or DRO problems. 
\subsection{Ambiguity Sets}
DRO can provide robust protection against scenarios generated by any distribution inside an ambiguity set $\mP$.
However, depending on the size of the set e.g., if  $\mP = \mathcal{P}_0$, the protection may be too conservative.
We can reduce this conservativeness by integrating the information gained from the data generated by the true distribution and thus shrinking the ambiguity set.
We construct the ambiguity sets as the intersection of the probability simplex and a data-driven set which contains the true distribution (with high confidence).
These data-driven sets can be prescribed by multi-dimensional confidence intervals or metrics such as the $\ell_2$-norm, kernel based deviation etc. 
These metrics are selected as they can provide guarantees about containing the true distribution. 
Their dependence on data also ensures that with more data, the sets
converge to the true distribution.

Our data-driven ambiguity sets~$\mP_{t}$, where $t\in \N$ equals the number of data points, can be of the following two forms
\begin{equation*}
	\mP_{t} \coloneqq \left\{p \in \mP \mid l_t\leq p \leq u_t \right\}\text{ or } \mP_{t}\coloneqq \left\{p \in \mP \mid d(p, \hat{p}_t) \leq \epsilon_t \right\}.
\end{equation*}
Here $l_t,u_t\in[0,1]^{|\mS|}$ are lower and upper bounds of confidence intervals, $\hat{p}_t$ is the empirical distribution estimator for $p^*$ and $d(\cdot, \cdot)$ and $\epsilon_t$ denote a distance metric and its respective bound. 
In both cases, the values of the parameters $l_t, u_t$ and $\epsilon_t$ are selected such that the true distribution lies inside the sets $\mP_t$ with high probability i.e.,
\begin{equation*}
\mathbb{P}\left(p^* \in \mP_t\right) \geq 1-\delta_t,
\end{equation*}
where $\delta_t\in (0,1).$
One key requirement for the ambiguity sets~$\mP_{t}$ over time~$t=0,...,T$ is that they contain
the true distribution inside \emph{all of them} with high probability, that means
\begin{align*}
p^* \in \bigcap_{t=0,...,T}\mathcal{P}_t,
\end{align*}
with a probability of at least $1-\delta$.
This is achieved by ensuring that each set $\mP_t$ in round $t$ contains the true distribution $p^*$ with probability at least $1 -\delta_t$ such that $\sum_{t=1}^{\infty}\delta_t < \infty$. 
In this paper, we choose $\delta_t = \frac{6 \delta}{\pi^2 t^2}$ for some predefined $\delta \in (0,1)$.
With this in mind, we have the following lemma.

\begin{lemma}\label{Lemma:Union_bound}
	For ambiguity sets constructed with confidence $\delta_t \coloneqq \frac{6 \delta}{\pi^2 t^2}$ and $\delta\in (0,1)$, it follows that the true data generating distribution  $p^* \in \bigcap_{t=0,...,T}\mathcal{P}_t$ with a probability of at least $1-\delta$.
\end{lemma}
\proof{Proof:}
Given a sequence of events $A_t$, we can estimate the probability of their intersection with Boole's inequality as follows, 
\begin{align*}
\mathbb{P}\left(\bigcap_{t=1,...,T} A_t\right) = 1 - \mathbb{P}\left(\bigcup_{t=1,...,T} A_t^c\right) \geq 1 - \sum_{t=1,...,T} \mathbb{P}( A_t^c).
\end{align*}
Let $A_t$ be the event that the true distribution $p^*$ lies inside the uncertainty set $\mathcal{P}_t$.
By the definition of an ambiguity set with confidence $\delta_t$, we know that $\mathbb{P}(A_t) \geq 1 - \delta_t$. This means that $\mathbb{P}(A_t^c) \leq \delta_t$.
This inequality and the limit of the over-harmonic series (2-series: $\sum_{t=1}^\infty \frac{1}{t^2}=\frac{\pi^2}{6}$) allows us to show that the probability of the event $p^* \in \bigcap_{t=1...T} \mathcal{P}_t$ (here without $t=0$) is at least
\begin{align*}
1- \sum_{t=1}^T \delta_t \geq 1 - \sum_{t=1}^\infty \delta_t = 1 - \frac{6 \delta}{\pi^2} \sum_{t=1}^\infty \frac{1}{t^2} =1- \delta.
\end{align*}
Since $$\bigcap_{t=1,...,T}\mathcal{P}_t \subset \mathcal{P}=\mathcal{P}_0 \Rightarrow \bigcap_{t=1,...,T}\mathcal{P}_t = \bigcap_{t=0,...,T}\mathcal{P}_t,$$ we conclude the proof (here with $t=0$). \Halmos

Note that the above result continues to hold as $T \ra \infty$. 
The above lemma shows that if the confidence probability $1-\delta_t$ increases fast enough, the inclusion of the true distribution $p^*$ in individual ambiguity sets $\mP_t$ is sufficient to guarantee uniform inclusion over all sets. 

\subsection{Choice of Ambiguity Sets}

The choice of ambiguity sets depends on the application, the historical information available and level of protection desired.
We discuss different ambiguity sets that are easily applicable and provide probabilistic performance guarantees.
These different ways of constructing data-driven ambiguity sets  demonstrate the flexibility of our approach.

\paragraph{Confidence Intervals.}
Confidence intervals for multinomial distributions can be calculated using various methods (see \citet{wang2008exact} for a survey).
We use the analytic formula from \citet{fitzpatrick1987quick} for the construction of multi-dimensional intervals $\mI_t\subseteq [0,1]^{|\mS|}$  via
\begin{equation}
	\mI_{kt}  = [l_{kt}, u_{kt}]\coloneqq  \left[\hat{p}_{kt} - \frac{z_{\frac{\delta_t}{2}}}{2 \sqrt t} \text{ , } \hat{p}_{kt} + \frac{z_{\frac{\delta_t}{2}}}{2 \sqrt t}\right], \label{Eq:Fitzpatrick}
\end{equation}
in each round $t=1,...,T$.
Here, $\hat{p}_t$ is the maximum likelihood estimator for $p^*$ and $z_{\frac{\delta_t}{2}}$ denotes the upper $(1 - \frac{\delta_t}{2})$-percentile of the standard normal distribution.  
The corresponding ambiguity sets can then be the intersection of the confidence intervals and the probability simplex, i.e.,
\begin{align*}
\mathcal{P}_t \coloneqq \mP_0 \cap \mI_t = \lrset{p \in \mathcal{P}_0 \mid l_t \leq p \leq u_t},
\quad t = 1, \ldots, T.
\label{label:pn}
\end{align*}
The parameters $l_t$ and $u_t$ can be updated as mentioned in the definition of $\mathcal{I}_{kt}$. 
The parameter $\hat{p}_t$ is the empirical probability distribution and can be updated by  counting the observations of each scenario up to $t$.
Confidence intervals work well in practice and are algorithmically preferable as they only impose linear constraints on the ambiguity sets.
This provides significant scalability to problems which use confidence intervals.
However, one disadvantage of them is that the probability guarantees that they provide only hold asymptotically. 
As a consequence, we extend our method to the following two sets which
provide finite sample guarantees while being easy to reformulate.

\paragraph{$\ell_2$-norm sets.}
For the finite sample setting, we have the following guarantee as proven in~\citet{weissman2003inequalities},
\begin{align*}
	\mathbb{P} \Big( \Vert \hat{p}_t - p^* \Vert_1 \leq \sqrt{\frac{2|\mS| \log 2/\delta_t}{t}}  \Big) \geq 1-\delta_t.
\end{align*} 
Given this inequality, along with the observation that $\|p\|_2 \leq \|p\|_1$ we can construct the following set
$$
\mP_t = \left\{p \in \mP \mid \|p - \hat{p}_t\|_2 \leq \epsilon_t \right\},
$$
with $\epsilon_t\coloneqq\sqrt{\frac{2|\mS| \log 2/\delta_t}{t}} $, which provides the containment guarantee $1-\delta_t$.
For the above ambiguity set, $\hat{p}_t$ and $\epsilon_t$ are the two parameters which define the set. 
These can be updated by counting the scenario observations and as per the definition of $\epsilon_t$ respectively. 

\paragraph{Kernel based ambiguity sets.}
Kernel based ambiguity sets are another alternative that provide similar finite sample probability guarantees while allowing for flexibility in how the different scenarios are weighted. Given a kernel function $k_M(s_i,s_j) : \mS\times \mS \rightarrow \R$ defined over the scenarios with kernel matrix $M$, we have the following probability guarantee~\citep{KirschnerBogunovicJegelkaKrause2020},
$$
\mathbb{P}\Big(\|\hat{p}_t - p^*\|_M \leq \frac{\sqrt{C}}{\sqrt{t}}(2 + \sqrt{2 \log(1/\delta_t)})\Big) \geq 1- \delta_t,
$$
if $k_M(s_i,s_j)\leq C$ for all scenarios. Here, $\|q\|_M := \sqrt{q^{\top}Mq}$ . 
With this inequality, we can construct an ambguity set similar to the $\ell_2$-norm case with $\epsilon_t := \frac{\sqrt{C}}{\sqrt{t}}(2 + \sqrt{2 \log(1/\delta_t)})$.
Like $\ell_2$ norm ambiguity sets, kernel based sets require two key parameters $\hat{p}_t$ and $\epsilon_t$. 
Both can be updated either by counting scenarios or the definition of $\epsilon_t$ for kernel sets. 
Using these ambiguity sets, the resulting~\eqref{Eq:DRON} forms a min-max problem, where the inner maximization problem optimizes a linear objective over a finite-dimensional convex feasible set.
Therefore, \eqref{Eq:DRON} can be equivalently reformulated using duality theory as is commonly applied in convex robust optimization. 
These reformulations are discussed in the following section.  

\subsection{Solving DRO via Reformulation}
When the ambiguity set is constructed with confidence intervals, the inner maximization problem forms a linear program.
Then, by strong duality, for confidence intervals, \eqref{Eq:DRON} is equivalent to
\begin{equation}
	\begin{aligned}
		\min_{x, z, \alpha, \beta} \text{ } & z
			- \langle l_t, \alpha\rangle + \langle u_t, \beta\rangle
			\\
		\st \text{ } & z - \alpha_k + \beta_k \geq f(x, s_k)
			\hspace{0.1cm}\hspace{0.2cm} \forall k = 1, \ldots, \abs{\mathcal{S}},
			\\
		& \alpha, \beta \geq 0, \\
		& x \in \mathcal{X}, \text{ } z \in \R,
			\text{ } \alpha, \beta \in \R^{\abs{\mathcal{S}}}.
	\end{aligned}
	\label{Eq:DRO_ref}%
\end{equation}
Here, the dual variables $ \alpha_k $ and $ \beta_k $
price the uncertainty for scenario $ s_k \in \mathcal S $.
This problem is
of the same problem class as~\eqref{Eq:Stoch_nominal},
however larger in size.
The reformulated DRO problem grows linearly with the number of
scenarios which may become prohibitive if the cardinality of
$\mathcal{S}$ is large.
Thus, the difficulty of solving~\eqref{Eq:DRO_ref} depends on the complexity of $f$ and the cardinality of~$\mathcal{S}$.

For the $\ell_2$-norm and kernel based ambiguity sets, the dual reformulation is given by
\begin{equation}
\label{eq:DRON_l2}
\begin{aligned}
	\hspace{-2mm}\min_{x,z, r} &\sum_{k=1}^{|S|} \hat{p}_{kt} \left(f(x,s_k)  + r_k \right)+ \epsilon_t \| \mathbf{f}(x) + z \mathbf{e}  + \mathbf{r} \|_{A^{-1}}\\
	\text{s.t.} &\; x \in \mathcal{X}, r_k \geq 0,
\end{aligned}
\end{equation}
where the matrix $A = I$ for the $\ell_2$-norm case and $M$ for the kernel based sets. 
The vectors $\mathbf{f}, \mathbf{e},$ and $\mathbf{r}$ denote the functions $f(x,s_k)\text{ for all } k$, the vector of all ones and the set of all $r_k$ respectively. 
Although norm constraints typically remain algorithmically
tractable, solving this reformulation is practically more difficult
than solving SO. 

\subsection{DRO over Time}
We now discuss our baseline framework of DRO over Time.
The iterative procedure outlined
in~Algorithm~\ref{Alg:DRO_over_time_draft} solves a sequence of
reformulated DRO problems while learning from data.
Prior to the first round, we assume that we have no data about the scenarios
and therefore initialize the ambiguity set
with the full probability simplex.
As more information comes in, the ambiguity set is updated. 

\begin{algorithm}
	\caption{DRO over Time}
	\begin{algorithmic}[1]
		\STATE \textbf{Input}: functions $ f(\cdot, \scenario) $
		for $ \scenario \in \scenarioset $, feasible set $ \iX $, initial ambiguity set $\mP_0$
		\STATE \textbf{Output}: sequence of DRO solutions $ x_1, \ldots, x_T $
		\FOR{$ t = 1 $ {\bfseries to} $T$}
		\STATE $ x_t \leftarrow $ solve Problem~\eqref{Eq:DRO_ref} or~\eqref{eq:DRON_l2} for $ \mathcal{P}_{t-1} $
		\STATE $ \mathcal{P}_{t} \leftarrow $ observe data
		and update set parameters such as $\hat{p}_t, l_t, u_t$ and $\epsilon_t$ as per the type of ambiguity set.
		\ENDFOR
	\end{algorithmic} \label{Alg:DRO_over_time_draft}
\end{algorithm}
The following theorems show the convergence of the DRO solution
$\widehat{J}_t$ to the true solution of the SO problem $J^*$.
The proofs are provided in the electronic companion as they are
adaptations from \cite{esfahani2018data} to our setting.
In the latter, 
analogous results have been proven for ambiguity sets
constructed with the Wasserstein metric. 
\begin{theorem}
	\label{thm:asymp_consis}
	If the feasible set $\mX$ is compact, then the optimal value of~\ref{Eq:DRON} converges over time to the optimal value of~\ref{Eq:Stoch_nominal} with probability 1, i.e.,
	\[\lim_{t \ra \infty} \widehat{J}_t = J^* \text{ with probability 1.}\]
\end{theorem}

\begin{theorem}
	\label{thm:soln_conv}
	Let $\{x_t\}_{t=1}^{\infty}$ be a sequence of  optimal solutions to the problem~\ref{Eq:DRON}. 
	If the feasible set $\mX$ is compact and the function $f(x,s)$ is continuous in $x$, then any accumulation point of $\{x_t\}_{t=1}^{\infty}$ is almost surely an optimal solution to the problem~\ref{Eq:Stoch_nominal}. 
\end{theorem}
These two results show the validity of the DRO over Time paradigm and guarantee that with a sufficient amount of data, Algorithm~\ref{Alg:DRO_over_time_draft} converges to the solution of the SO problem with the true distribution.
This is important as this indicates the importance of gathering more data to obtain better solutions.

In Algorithm~\ref{Alg:DRO_over_time_draft},
the DRO problem has to be solved repeatedly
in each round.
This is not viable for large problems as it requires long computation times. 
We remedy this by introducing a new algorithm
which approximates the DRO problem over time and updates the ambiguity sets via online learning.
For this online algorithm, we can show a dynamic regret bound that bounds the error of the approximation with sublinear expression in the number of rounds.


%% file: dro_via_online_learning.tex
\section{Online Robust Optimization}
\label{Sec:DRO_over_time}
\label{sec:drotime}
In this section, we introduce the online learning and optimization algorithm which is main contribution of our work. 
As in~\cite{ben-tal:oracle_ro} and~\cite{pokutta2021adversaries}, we consider robust optimization as a game between two players. 
The online algorithm can then be roughly described as alternating between solving the optimization problem for each player given the solution of the other.
Thus, for each round $t=1,\dots,T$, we decompose the min-max problem into two subproblems (one for each player) and perform the following steps:
\begin{enumerate}
	\item First, the \emph{$p$-player} determines $p_{t}$ via an appropriate  algorithm applied to problem
	$$
	\max_{p\in \mathcal{P}_{t-1}} \mathbb{E}_{s \thicksim p}\left[f(x_{t-1},s)\right]\label{Eq:Adversarial},
	$$ 
	based on the solution $x_{t-1}$ from the previous round.
	\item Then, the \emph{$x$-player} computes  $x_{t}$ as a solution using the previously calculated $p_t$
	$$
	\min_{x\in \mathcal{X}} \text{ }  \mathbb{E}_{s \thicksim p_{t}}\left[f(x,s)\right]. \label{Eq:Stoch}
	$$	
\end{enumerate}
Due to the fact that we consider probability distributions over a finite scenario set, the optimization problem of the $p$-player is finite-dimensional and therefore the Online Gradient Descent~\citep{zinkevich2003online} is a canonical choice here as a learning algorithm. 
Given $p_{t-1}$ and $x_{t-1}$, 
the update rule consists of a descent step 
\begin{align*}
\tilde{p}_{t} = p_{t-1} + \eta \nabla_p \mathbb{E}_{s \sim p_{t-1}}\left[f(x_{t-1},s)\right],
\end{align*}
with step size $\eta > 0$ and a subsequent projection step to ensure feasibility
\begin{align*}
p_{t} = \arg\min_{p\in \mathcal{P}_{t-1}} \textstyle\frac{1}{2}\Vert p - \tilde{p}_{t} \Vert^2.
\end{align*}
The probability distribution in the next iteration is therefore given as the unique solution
\begin{align*}
p_{t} = \arg\min_{p\in \mathcal{P}_{t-1}}& \text{ }\left\langle -\eta \nabla_p \mathbb{E}_{s \thicksim p_{t-1}}\left[f(x_{t-1},s)\right], p \right\rangle 
 + \frac{1}{2}\Vert p - p_{t-1} \Vert ^2.
\end{align*}
Note that the $x$-player uses $p_t$ in round $t$  to estimate $x_t$ and solves a standard SO problem which is easier to solve than the reformulated DRO program. 
This $p_t$ is then used to estimate $p_{t+1}$ in round $t+1$.
As such, the value of $x_t$ depends on $p_t$ in round $t$ (thus
$x_t$ and $p_t$ are not conditionally independent in round $t$). 
Therefore, as per~\cite{pokutta2021adversaries}, the learner for $p_t$ has to be a strong learner in order to ensure sublinear regret. 
This paper also shows that the online gradient descent algorithm satisfies the necessary conditions for it to be a strong learner.

\subsection{Algorithm}
We provide a pseudo code of our method for DRO over Time via online robust
optimization in Algorithm~\ref{Alg:DRO_via_adversarial}.
It combines alternating solutions of the $\min$ and $\max$ problems with the update of the ambiguity sets.
For the sake of simplicity, we assume that our algorithm starts without any knowledge of the probability distribution over the scenarios and therefore we initialize the ambiguity set as the full probability simplex in step~1.
The algorithm can be easily modified to incorporate any historical information.
The initialization of $p_0\in\mP_{0}$ and $x_0\in\mX$ in step~4 and step~5 can be chosen arbitrarily and does not effect our theoretical results.
Each round $t=1,..., T$ starts with the update of $p_t$ via projected
gradient descent and of $x_t$ as the solution of an SO
problem.
At the end of the round, we observe new data in form of (i.i.d.) scenario observations and update the
ambiguity set as explained in Section~\ref{sec:dddro}.
\begin{algorithm}[htb]
	\caption{DRO over Time with Online Projected Gradient Descent}
	\begin{algorithmic}[1]
		\STATE {\bfseries Input:} functions $ f(\cdot, s) $ for $ s \in \iS $, feasible set $\iX$, initial ambiguity set $\mP_0$
		\STATE {\bfseries Output:}  $x_{1},\ldots,x_{T}$
		\STATE Set $\mathcal{P}_0 = 	\set{p \in [0, 1]^{\abs{\mathcal{S}}} \mid
			\sum_{k = 1}^{\abs{\mathcal{S}}} p_k = 1} $
		\STATE Set $p_0=\left(\frac{1}{|\mathcal{S}|}, ...,\frac{1}{|\mathcal{S}|} \right)\in \left[0,1\right]^{|\mathcal{S}|}$
		\STATE Set $x_0=\min_{x\in \mathcal{X}} f(x,s_1)  $
		\FOR{$t=1$ {\bfseries to} $T$}
		\STATE $\tilde{p}_{t} \leftarrow p_{t-1} + \eta \nabla_p \mathbb{E}_{s \sim p_{t-1}}\left[f(x_{t-1},s)\right]$
		\STATE 	$p_{t} \leftarrow \arg\min_{p\in \mathcal{P}_{t-1}} \textstyle\frac{1}{2}\Vert p - \tilde{p}_{t} \Vert^2$ 
		\STATE $x_{t} \leftarrow \arg\min_{x\in \mathcal{X}} \mathbb{E}_{s \sim p_{t}}\left[f(x,s)\right] $
		\STATE $\mathcal{P}_{t} \leftarrow$ observe data
		and update set parameters such as $\hat{p}_t, l_t, u_t$ and $\epsilon_t$ as per the type of ambiguity set.
		\ENDFOR
	\end{algorithmic}\label{Alg:DRO_via_adversarial}
\end{algorithm}

Algorithm~\ref{Alg:DRO_via_adversarial} provides a sequence of solutions $x_t$ for each time step $t=0,...,T$.
In order to prove that the quality of the solutions $x_t$ improves over time, we bound the average gap over time between the worst case performance of $x_t$ and the optimal worst case (exact DRO) solution. 
Since the feasible set of the $p$-player changes over time it is not possible to apply existing techniques for regret bounds on min-max problems to our setting as these techniques have primarily focused on stationary ambiguity sets. 
In this paper, we extend these existing techniques to the case of shrinking ambiguity sets by leveraging the fact that all the ambiguity sets contain a common distribution (the true distribution). 


The two main ingredients for the theoretical analysis are a constrained gradient ($\nabla_p \mathbb{E}_{s \sim p}\left[f(x,s)\right]$) and a constrained path length  ($\sum_{t=1}^{T}\frac{1}{2}\Vert p_t - q_t \Vert^2\leq h(T)$ for all $p_t\in\mP_{t-1}, q_t\in\mP_{t}$) for the online gradient descent to work on non-stationary feasible sets.   
The former is a classical assumption for steepest descent algorithms while the latter is commonly used in dynamic regret bounds for online algorithms, see e.g.~\citep{path_length_dynamic_regret, dynamic_regret}.


For constant ambiguity sets, it is known that for the static regret (which compares the performance difference of online solutions to a single best action in hindsight) bounds of the  
form $\mathcal{O}({1}/{\sqrt{T}})$ can be derived~\citep{pokutta2021adversaries, besbes2015non}. 
For our case, a careful analysis of the algorithm leads to a dynamic regret bound of $\mathcal{O}({\sqrt{h(T)}}/{\sqrt{T}})$ that is presented in
Theorem~\ref{Th:Regret_bound}, with
the corresponding bounding terms $h(T)$ for the path lengths being proven afterwards.
This is consistent with other findings about dynamic regret bounds in the literature on online learning, see~\citep{dynamic_regret}.
We are able to achieve a bound for this new setting with shrinking ambiguity sets.

\begin{theorem}[Dynamic regret bound] \label{Th:Regret_bound}
	Let $f:\mathcal{X}\times\mathcal{S}\rightarrow\mathbb{R}$ be uniformly bounded, i.e., for all $(x,s) \in \mathcal{X}\times\mathcal{S}$, there exists a constant $G>0$ such that $|f(x,s)|\leq G$.  Let $\eta := \sqrt{\frac{2 h(T)}{G^2 T|\mS|}}$ where  $\sum_{t=1}^{T}\frac{1}{2}\Vert p_t - q_t \Vert^2 \leq  h(T)$ for $p_t \in \mP_{t-1}$ and $q_t \in \mP_t$. 
	The output $(x_1,...,x_{T})$ from Algorithm~\ref{Alg:DRO_via_adversarial} with confidence update $\delta_t \coloneqq \frac{6 \delta}{\pi^2 t^2}$ and $\delta\in (0,1)$ fulfills
	\begin{align*}
	\frac{1}{T} &\sum_{t=1}^{T}  \left(\max_{p \in \mP_{t}}\mathbb{E}_{s \thicksim p}\left[f(x_{t},s)\right]  - \min_{x\in \mathcal{X}} \max_{p\in \mathcal{P}_{t}} \mathbb{E}_{s \thicksim p}\left[f(x,s)\right] \right) 
	 \leq G \sqrt{\frac{2|\mS|h(T)}{T}} + \frac{2G}{T},
	\end{align*}
	with probability at least $1-\delta$.
\end{theorem}
\proof{Proof of Theorem \ref{Th:Regret_bound}}
Define $g_t(p)\coloneqq -\mathbb{E}_{s \thicksim p}\left[f(x_t,s)\right]$. An online gradient descent iteration is given by
\begin{align*}
	p_{t+1} = \arg\min_{p\in \mathcal{P}_t} \text{ }\left\langle \eta \nabla g_t(p_t), p \right\rangle + \frac{1}{2}\Vert p - p_{t} \Vert^2,
\end{align*}
with the variational inequality
$$
\left\langle \eta \nabla g_t(p_t), u_t - p_{t+1}  \right\rangle + \left\langle p_{t+1} - p_{t}, u_t - p_{t+1}  \right\rangle \geq 0, \text{ for all  }u_t \in \mathcal{P}_{t}
$$
as the optimality criteria.
Classical theory 
for gradient descent (by rearranging the previous inequality and using Cauchy-Schwarz) yields
\begin{align*}
	\left\langle \eta \nabla g_t(p_t), p_t - u_t  \right\rangle \leq &\frac{1}{2}\Vert p_{t} - u_t \Vert^2 - \frac{1}{2}\Vert p_{t+1}-u_t  \Vert^2 
	+ \frac{\eta^2}{2} \Vert \nabla g_t(p_t) \Vert ^2.
\end{align*}
Summation over rounds $t=1,...,T$ results in the following inequality	for all $u_t \in \mathcal{P}_{t}$:
\begin{align*}
	&\sum_{t=1}^T\left\langle \eta \nabla g_t(p_t), p_t - u_t  \right\rangle \leq  \sum_{t=1}^T\frac{\eta^2}{2} \Vert \nabla g_t(p_t) \Vert ^2 
	+ 	\sum_{t=1}^T \left(\frac{1}{2}\Vert p_{t} - u_t \Vert^2 - \frac{1}{2}\Vert p_{t+1} - u_t \Vert^2 \right). 
\end{align*}
Next, we bound the terms on the right side starting with 
\begin{align*}
	& \sum_{t=1}^T \left(\frac{1}{2}\Vert p_t - u_t \Vert^2 - \frac{1}{2}\Vert p_{t+1} - u_t \Vert^2 \right)
	\leq \sum_{t=1}^T \frac{1}{2}\Vert p_t - u_t \Vert^2.
\end{align*}
Note that this is, in contrast to classical steepest descent theory, not a telescoping sum. 
We make use of the fact that $\sum_{t=1}^{T}\frac{1}{2}\Vert p_t - u_t \Vert^2 \leq  h(T)$ for all $u_t \in \mathcal{P}_{t}$ and all rounds $t=1,...,T$ with probability at least $1-\delta$ as per Lemma~\ref{lemma:all_path_lengths}. 
We provide the proofs for these path length bounds in the next extra subsection.
Furthermore, since $g_t(p)=-\mathbb{E}_{s\thicksim p}\left[f(x_{t},s)\right]$ is linear in $p$, we use the gradient bound
\begin{align*}
	\Vert \nabla g_t(p_t) \Vert^2 = \sum_{k = 1}^{|\mathcal{S}|} |f(x_t,s_k)|^2 \leq |\mS|G^2,
\end{align*}
for all $t = 1, ..., T$, because $f$ is bounded.
Thus we get 
\begin{align*}
	\sum_{t=1}^T&\left\langle \nabla g_t(p_t), p_t - u_t  \right\rangle \leq \frac{ h(T)}{\eta} + \frac{\eta}{2} T |\mS| G^2, 
\end{align*}
for all $u_t \in \mathcal{P}_{t},$ for all $t = 1, ..., T$ with probability at least $1-\delta$.
Choosing the bound minimizing step size (minimize right-hand side with respect to $\eta$) 
$\eta:=\sqrt{\frac{2 h(T)}{G^2 |\mS| T}}$ yields
\begin{align*}
	\sum_{t=1}^T\left\langle \nabla g_t(p_t), p_t - u_t  \right\rangle \leq G \sqrt{2|\mS|h(T) T}.
\end{align*} 
Since $g_t(p)=-\mathbb{E}_{s\thicksim p}\left[f(x_{t},s)\right]$ is linear in $p$ for all $t=1,...,T$, it follows
\begin{align*}
	&\sum_{t=1}^T \left( \mathbb{E}_{s\thicksim u_t}\left[f(x_{t},s)\right] - \mathbb{E}_{s\thicksim p_t} \left[f(x_t, s)\right] \right) 
	= \sum_{t=1}^T \left\langle \nabla g_t(p_t), p_t - u_t \right\rangle \leq G \sqrt{2|\mS| h(T) T}.
\end{align*}
Now we choose in each round $t=1,...,T$ the worst-case $ 
u_t \coloneqq \arg \max_{p \in \mP_{t}} \mathbb{E}_{s\thicksim p}\left[f(x_{t},s)\right] \in \mP_t $
and recall the definition of $x_{t} = \arg\min_{x\in \mathcal{X}} \mathbb{E}_{s \sim p_{t}}\left[f(x,s)\right]$ to obtain 
\begin{align*}
	\sum_{t=1}^T &\left( \mathbb{E}_{s\thicksim u_t}\left[f(x_t,s)\right] - \mathbb{E}_{s\thicksim p_t} \left[f(x_t, s)\right] \right) 
	= \sum_{t=1}^T \left( \max_{p \in \mP_t}\mathbb{E}_{s\thicksim p}\left[f(x_t,s)\right] - \min_{x \in \mX} \mathbb{E}_{s\thicksim p_t} \left[f(x, s)\right] \right).
\end{align*}
Since $p_t \in \mP_{t-1}$, we know that $\min_{x \in \mX} \mathbb{E}_{s\thicksim p_t} \left[f(x, s)\right] \leq \min_{x \in \mX} \max_{p \in \mP_{t-1}}\mathbb{E}_{s\thicksim p} \left[f(x, s)\right]$ for all $t=1,...,T$  and thus we can conclude
\begin{align*}
	\sum_{t=1}^T& \left( \max_{p \in \mP_t} \mathbb{E}_{s\thicksim p}\left[f(x_t,s)\right] - \min_{x \in \mX} \max_{p \in \mP_{t-1}} \mathbb{E}_{s\thicksim p} \left[f(x, s)\right] \right)
	\leq G \sqrt{{2|\mS| h(T) T}},
\end{align*}
with a probability of at least $1-\delta$.
We add and subtract $\min_{x \in \mX} \max_{p \in \mP_t} \mathbb{E}_{s\thicksim p} \left[f(x, s)\right]$ on the LHS.  Rearranging the terms like this allows us to write the LHS as
\begin{align*}
	&\sum_{t=1}^T \left( \max_{p \in \mP_t} \mathbb{E}_{s\thicksim p}\left[f(x_t,s)\right] - \min_{x \in \mX} \max_{p \in \mP_t} \mathbb{E}_{s\thicksim p} \left[f(x, s)\right] \right)\\ & + \sum_{t=1}^T \min_{x \in \mX} \max_{p \in \mP_t} \mathbb{E}_{s\thicksim p} \left[f(x, s)\right]  -  \min_{x \in \mX} \max_{p \in \mP_{t-1}} \mathbb{E}_{s\thicksim p} \left[f(x, s)\right]. 
\end{align*}
The last two terms telescope. Bringing them to the RHS and using the upper bound $G$ on $|f(x,s)|$ we can conclude
\begin{align*}
	\sum_{t=1}^T& \left( \max_{p \in \mP_t} \mathbb{E}_{s\thicksim p}\left[f(x_t,s)\right] - \min_{x \in \mX} \max_{p \in \mP_t} \mathbb{E}_{s\thicksim p} \left[f(x, s)\right] \right)
	\leq G \sqrt{2|\mS| h(T) T} + 2G \;\;\;\;\;\text{ w.p. } 1 - \delta.
\end{align*}
Dividing by $T$ on both sides completes the proof. \Halmos
A second regret bound is also proven in the electronic companion which allows for better dependence on the number of scenarios $|\mS|$  but at a cost of worse dependence on the total iteration count $T$.
This is achieved by replacing the path length $\sum_{t=1}^T \frac{1}{2}\Vert p_{t} - u_t \Vert^2\leq h(T)$ by $ \sum_{t=2}^{T}\|p_t - q_t\| \leq h'(T) $, which leads to different bounds. 
Numerically, it is better in the beginning i.e., for small $t$ but has worse asymptotic behavior. 

As stated in the following Corollary, we also observe that for all considered ambiguity sets the dynamic regret bound converges to zero.
\begin{corollary}[Convergence of Regret]
	\label{thm:regret_conv}
	If $\lim_{T \ra \infty}h(T)/T = 0$, the dynamic regret converges to 0 with probability $1-\delta$ i.e.,
	\begin{align*}
	\lim_{T \ra \infty} \frac{1}{T} \sum_{t=1}^{T}& \Big( \max_{p \in \mP_{t}}\mathbb{E}_{s \sim p}\left[f(x_{t},s)\right] 
	- \min_{x\in \mathcal{X}} \max_{p\in \mathcal{P}_{t}} \mathbb{E}_{s \sim p}\left[f(x,s)\right] \Big) = 0.
	\end{align*}
\end{corollary}

\proof{Proof of Corollary~\ref{thm:regret_conv}}
We know that 
\begin{align*}
\limsup_{T \ra \infty} \frac{1}{T} \sum_{t=1}^{T} &\left( \max_{p \in \mP_{t}}\mathbb{E}_{s \sim p}\left[f(x_{t},s)\right] \right. - \left. \min_{x\in \mathcal{X}} \max_{p\in \mathcal{P}_{t}} \mathbb{E}_{s \sim p}\left[f(x,s)\right] \right) \\
\leq \limsup_{T \ra \infty} & \left(G \sqrt{\frac{2h(T)|\mS|}{T}} + \frac{2G}{T}\right)= 0,
\end{align*}
with probability at least $1-\delta$. Thus, we can write,
\begin{align}
\label{eq:dyrub}
\limsup_{T \ra \infty} \frac{1}{T} \sum_{t=1}^{T}& \left( \max_{p \in \mP_{t}}\mathbb{E}_{s \sim p}\left[f(x_{t},s)\right] \right. - \left. \min_{x\in \mathcal{X}} \max_{p\in \mathcal{P}_{t}} \mathbb{E}_{s \sim p}\left[f(x,s)\right] \right) \leq  0.
\end{align}
To prove the lower bound,
let $\bar{x}_{t}$ be the optimal solution to the problem $\min_{x\in \mathcal{X}} \max_{p\in \mathcal{P}_{t}} \mathbb{E}_{s \sim p}\left[f(x,s)\right]$.
Then we can write the inner term in the left hand side (LHS) in equation~\eqref{eq:dyrub} as 
\begin{align*}
\max_{p \in \mP_{t}}\mathbb{E}_{s \sim p}\left[f(x_{t},s)\right] - \max_{p\in \mathcal{P}_{t}} \mathbb{E}_{s \sim p}\left[f(\bar{x}_{t},s)\right]. 
\end{align*}
We know that $\bar{x}_{t}$ is the optimal solution to the problem $\min_{x\in \mathcal{X}} \max_{p\in \mathcal{P}_t} \mathbb{E}_{s \sim p}\left[f(x,s)\right]$, this means that 
\begin{align*}
\max_{p\in \mathcal{P}_{t}} \mathbb{E}_{s \sim p}\left[f(x_t,s)\right] &\geq \min_{x\in \mathcal{X}} \max_{p\in \mathcal{P}_{t}} \mathbb{E}_{s \sim p}\left[f(x,s)\right] \\&= \max_{p\in \mathcal{P}_{t}} \mathbb{E}_{s \sim p}\left[f(\bar{x}_{t},s)\right]. 
\end{align*}
Thus we get 
\begin{align}
\label{eq:zero_lower_bound}
\max_{p\in \mathcal{P}_{t}} \mathbb{E}_{s \sim p}\left[f(x_t,s)\right] - \max_{p\in \mathcal{P}_{t}} \mathbb{E}_{s \sim p}\left[f(\bar{x}_{t},s)\right] \geq 0.
\end{align}
The above lower bound and the $ \limsup $-bound~\eqref{eq:dyrub}
together prove the result.  \Halmos

From the above result 
and Theorem~\ref{Th:Regret_bound}, we can observe that the dynamic regret, i.e., the average gap between the best solution in hindsight of each round and the solution evaluated in the algorithm, decreases at a rate of $\mathcal{O}(\sqrt{h(T)}/\sqrt{T})$ and tends to zero.
Therefore, we have a performance guarantee (with sublinear regret) when solving the DRO problems approximately over time. 
At the same time, Algorithm~\ref{Alg:DRO_via_adversarial} is applicable to large-sized problems in contrast to the reformulated DRO problems~\eqref{Eq:DRO_ref} or~\eqref{eq:DRON_l2}.

\subsection{Bounded Path Lengths}
The following lemma illustrates the high-probability path lengths for the confidence intervals, the kernel-based and $l_2$-norm ambiguity sets. 
\begin{lemma}\label{lemma:all_path_lengths}
	Given ambiguity sets of the form specified in Section~\ref{sec:dddro}, we have 
	\begin{align*}
	\frac{1}{2}\sum_{t=1}^{T}\Vert p_t - q_t \Vert^2 \leq h(T),
	\end{align*}
	for all $p_t\in\mP_{t-1}, q_t \in \mP_t$ with probability at least $1-\delta$.
	The functions $h(T)$ for different categories of ambiguity sets are as given: 
\begin{enumerate}
	\item \textbf{Confidence Intervals}: 
	$$h(T) = 8 |\mS| \log(\pi T) ({2} + \log T).$$
	\item \textbf{Kernel based ambiguity sets}:
	$$h(T) =  {\frac{1}{2}\left(2+ \frac{4\sqrt{C}}{\lambda}\right)^2}+ \frac{32C}{\lambda^2 } \log\frac{\pi T}{\sqrt{6\delta}}(1 + \log T),$$
	where $\lambda$ denotes the smallest eigenvalue of the kernel matrix $M$.
	\item \textbf{$\ell_2$-norm ambiguity sets:}
	$$
	h(T) = 8 |\mS|\log \frac{\pi T}{\sqrt{3 \delta}} ({2} + \log T).
	$$
\end{enumerate}
\end{lemma}
\proof{Proof of Lemma~\ref{lemma:all_path_lengths}.} 

\textbf{Confidence Intervals}.
We show in the electronic companion (Lemma~\ref{lemma:shrinking_difference}) that the ambiguity sets $\mP_t$ derived from \eqref{Eq:Fitzpatrick} with confidence update $\delta_t \coloneqq \frac{6 \delta}{\pi^2 t^2}$ and $\delta\in (0,1)$ for all rounds $t=1,...,T$ fulfill 
	\begin{align*}
	\sup_{x\in \mP_{0}, y\in \mP_{1}}\Vert x - y \Vert \leq \sqrt{16|\mS| \log\pi}
	\text{ and }
	\sup_{x\in \mP_{t-1}, y\in \mP_{t}}\Vert x - y \Vert \leq \frac{ \sqrt{16|\mS| \log(\pi {(t-1)})}}{\sqrt{{t-1}}},
\end{align*}
with a probability of at least $1-\delta$. This allows for calculating the function $h(T)$: 
{
\begin{align*}
	\frac{1}{2} \sum_{t=1}^{T}\|p_t - q_t\|^2 &\leq  \frac{1}{2} 16 |\mS| \log\pi + \frac{1}{2} \sum_{t=2}^{T}16 |\mS| \frac{\log(\pi (t-1))}{t-1}\\
	&\leq 8 |\mS| \log\pi + 8 |\mS| \log(\pi (T-1)) \sum_{t=1}^{T-1}\frac{1}{t}\\
	&\leq 8 |\mS| \log\pi + 8 |\mS| \log(\pi (T-1)) (1 + \log (T-1))\\
	&\leq 8 |\mS| \log(\pi T) (2 + \log T).
\end{align*}
}
The second and third inequalities are from bounding $t$ and from observing that $\sum_{t=1}^{T-1}(1/t) \leq 1 + \log(T-1)$. 

\noindent \textbf{Kernel based ambiguity sets}.
We show in the electronic companion (Lemma~\ref{lemma:kernel_shrinking_difference}) that given an ambiguity set of the form $\mP_t = \left\{p \in \mP \mid \|p - \hat{p}\|_M \leq \epsilon_t \right\}$ with $\epsilon_t\coloneqq \frac{\sqrt{C}}{\sqrt{t}}(2 + \sqrt{2 \log(1/\delta_t)}) $ with $\delta_t = \frac{6 \delta}{\pi^2 t^2}$ 
we have for $t \geq 2$, 
\begin{align*}
&\sup_{x \in \mP_{t-1}, y \in \mP_t} \|x - y\|_2 \leq \frac{8\sqrt{C}}{\lambda \sqrt{t-1}} \sqrt{ \log(\pi t/\sqrt{6\delta})},\\
&\text{ and } \sup_{x \in \mP_0, y \in \mP_1} \|x - y\|_2 \leq 2+\frac{4\sqrt{C}}{\lambda}\;\; \text{ for } t = 1, 
\end{align*}
with probability at least $1-\delta$.
Calculating the function $h(T)$, we have

\begin{align*}
	\frac{1}{2} \sum_{t=1}^{T}\|p_t - q_t\|_2^2 &\leq {\frac{1}{2}\left(2+ \frac{4\sqrt{C}}{\lambda}\right)^2}+  \sum_{t=2}^{T}\frac{32 C}{\lambda^2 ({t-1})}\log(\pi {t}/\sqrt{6\delta})\\
	&\hspace{-10mm}\leq {\frac{1}{2}\left(2+ \frac{4\sqrt{C}}{\lambda}\right)^2}+ \frac{32C}{\lambda^2 } \log(\pi T/\sqrt{6\delta})\sum_{t=2}^{T}\frac{1}{t-1}\\
	&\hspace{-10mm}\leq {\frac{1}{2}\left(2+ \frac{4\sqrt{C}}{\lambda}\right)^2}+ \frac{32C}{\lambda^2 } \log(\pi T/\sqrt{6\delta})(1 + \log T).
\end{align*}
Here, the first inequality arises from Lemma~\ref{lemma:kernel_shrinking_difference}. The second and third inequalities are from bounding $t$ and from observing that $\sum_{t=1}^{T-1}(1/t) \leq 1 + \log(T-1)$. 

\noindent \textbf{$\ell_2$-norm ambiguity sets}.
We show in the electronic companion (Lemma~\ref{lemma:l2_shrinking_difference}) that given an ambiguity set of the from $\mP_t = \left\{p \in \mP \mid \|p - \hat{p}\|_2 \leq \epsilon_t \right\}$ with $\epsilon_t\coloneqq\sqrt{\frac{2|\mS| \log (2/\delta_t)}{t}} $ and $\delta_t = \frac{6 \delta}{\pi^2 t^2}$, we have
$$\sup_{x \in \mP_{0}, y \in \mP_1} \|x - y\|_2 \leq 4\sqrt{|\mS| \log (\pi / \sqrt{3 \delta})}
\text{ and }
\sup_{x \in \mP_{t-1}, y \in \mP_t} \|x - y\|_2 \leq 4\sqrt{\frac{|\mS| \log (\pi (t-1) / \sqrt{3 \delta})}{{t-1}}},$$
with probability at least $ 1-\delta$.
For bound $h(T)$, it follows
\begin{align*}
	\frac{1}{2} \sum_{t=1}^{T}\|p_t - q_t\|^2 &\leq \frac{1}{2}(4\sqrt{|\mS| \log (\pi / \sqrt{3 \delta})})^2 + \frac{1}{2} \sum_{t=2}^{T}16 |\mS| \frac{\log(\pi ({t-1})/\sqrt{3 \delta})}{t-1}\\
	&\leq 8 |\mS| \log (\pi / \sqrt{3 \delta}) + 8 |\mS| \log(\pi T/\sqrt{3\delta}) \sum_{t=2}^{T}\frac{1}{t-1}\\
	&\leq 8 |\mS| \log(\pi T/\sqrt{3\delta}) (2 + \log T).
\end{align*}
  \Halmos
From Lemma~\ref{lemma:all_path_lengths}, we know that the shrinking ambiguity sets yield path length bounds $h(T)$ that increase at a rate of $\mathcal{O}(\log^2 T)$.
Therefore, the dynamic regret bound converges to zero with $\mathcal{O}(\log T / \sqrt{T})$.
It is only slightly worse than the typical $\mathcal{O}({1}/{\sqrt{T}})$ bounds.
In summary, this novel algorithm yields the remarkable benefit that
it integrates changing ambiguity sets (which represent growing
knowledge of the uncertainty) into regret bounds for online robust optimization.


%% file: numerical_results.tex
\section{Numerical Results}
\label{sec:num_res}

We illustrate the performance 
of Algorithm~\ref{Alg:DRO_via_adversarial} through numerical
experiments on mixed-integer linear and quadratic programs (MIPs \&
MIQPs) as well as on two real-world applications, namely
distributionally robust network design \citep{network_design} and
optimal route choice. 
These problems allow to demonstrate the wide applicability of our
novel approach.
The latter does not assume anything about the
structure of the original problem apart from the fact that it should
be algorithmically tractable, with an available solution approach.
We are thus able to apply DRO over time to both discrete and continuous optimization problems to account for various application types. 


All computations are carried out using a Python 3.8.5 implementation on  machines with Intel Core i7 CPU 2.80
GHz processor and 16 GB RAM. Each of these has four cores of 3.5 GHz each and 32 GB RAM.
We utilized SCIP 7.0 as the MIP and
MIQP solver \citep{scip} and IPOPT 3.13.3  \citep{ipopt} to compute the projections in Algorithm \ref{Alg:DRO_via_adversarial} as the solution of a convex optimization problem.
For the ambiguity sets, we choose $1-\delta=0.9$.
Different values for this parameter impact the size of the resulting
sets, but have no significant impact on the ability of Algorithm~\ref{Alg:DRO_via_adversarial} to learn DRO solutions. 
For the kernel-based ambiguity sets, we use the Gaussian kernel function $k_M(s_i,s_j) = \exp\left({-\frac {\Vert s_i - s_j \Vert_2^2} 2}\right)$.

\subsection{Benchmark Instances}\label{Sec:benchmark_instances}

To validate the performance of the method, we use publicly available
instances from well known benchmark libraries MIPLIB \citep{miplib2017} and QPLIB \citep{qplib}, which contain a collection of MIPs and MIQPs, respectively.
In all test cases, the uncertain objective function $f$ has the form
$f(x,s) = x^\top Qx + (c+s)^\top x + d$,
with $Q\in \R^{n\times n}$, $c\in \R^n$, $d\in \R$, $x\in
\mathcal{X}$ and $s\in \mathcal{S}$. 
For linear problems, we have $Q=0$. 
We generate different cost scenarios $s\in \mathcal{S}\subset \R^n, n>0$ by perturbing the coefficients of variables in the objective function randomly by up to 50\%.

We sort the
problems in the libraries by increasing number of variables and choose the first 15 MIPs and the first 10 convex MIQPs with bounded variables that can be solved within an hour of CPU time.
All instances are listed in the electronic companion.
The number of scenarios in the linear test cases is set to $|\mathcal{S}| \in \{10, 50\}$.
As the MIQP instances are more difficult to solve, we set
$|\mathcal{S}|=2$.

We numerically evaluate the dynamic regret
from Theorem \ref{Th:Regret_bound} and the worst-case expected value of solutions from Algorithm~\ref{Alg:DRO_via_adversarial} in comparison with the DRO solution obtained via reformulation \eqref{Eq:DRO_ref}.
In each round, one realization of the uncertain parameter is revealed.
These realizations are drawn randomly according to the initially generated probability distribution.
Thus, the number of data points equals the number of rounds.

\begin{figure}
	\centering
	\subfigure[Evolution of worst-case expected objective value for all three types of ambiguity sets in each round. In addition to the actual values, the worst-case average is given as a line.]{
		\includegraphics[width=8cm]{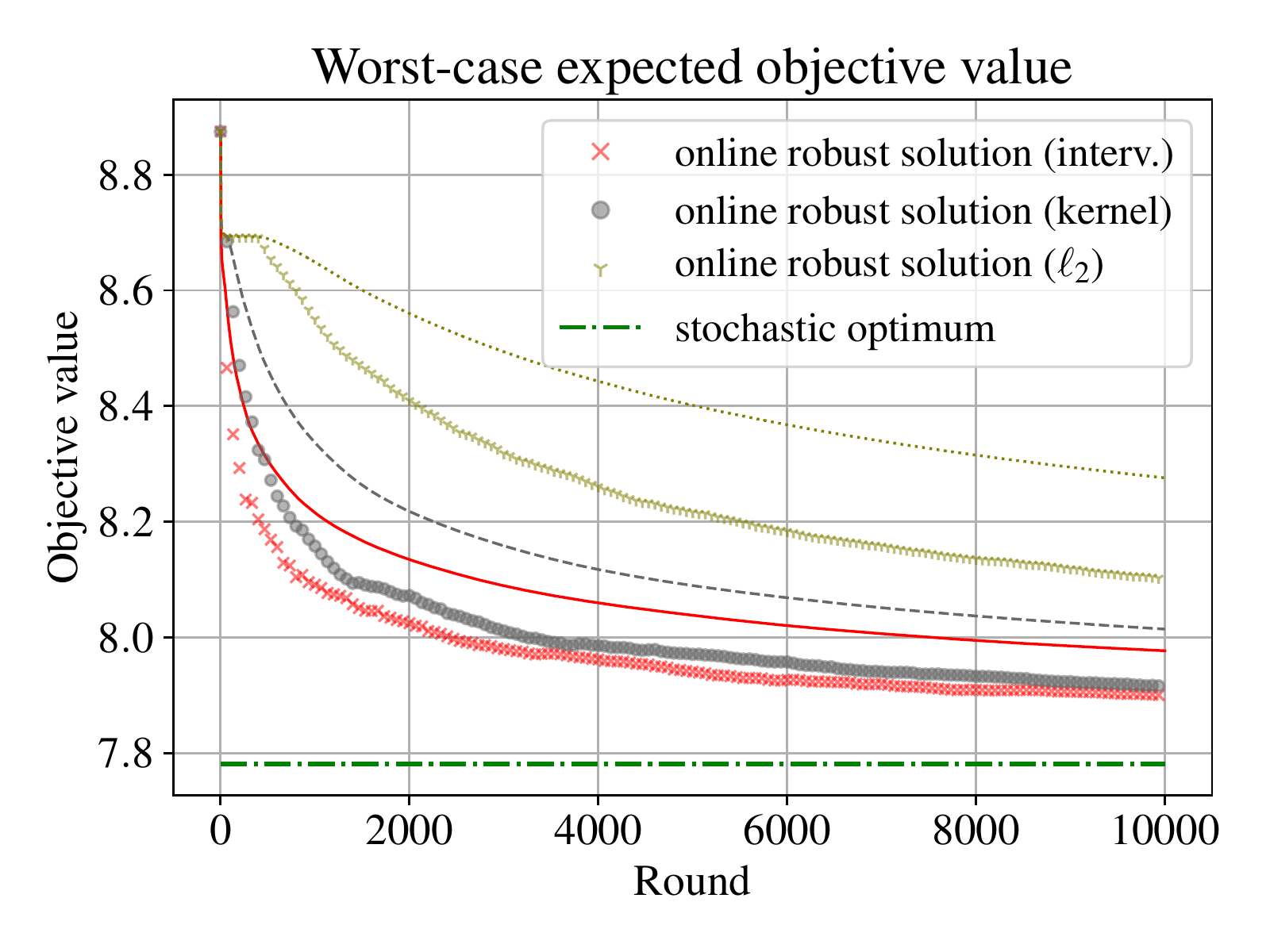}
		\label{fig:mip_regret1}	
	}
	\hspace{0.1mm}
	\subfigure[Convergence of online average to DRO average (dynamic regret). ]{\includegraphics[width=8cm]{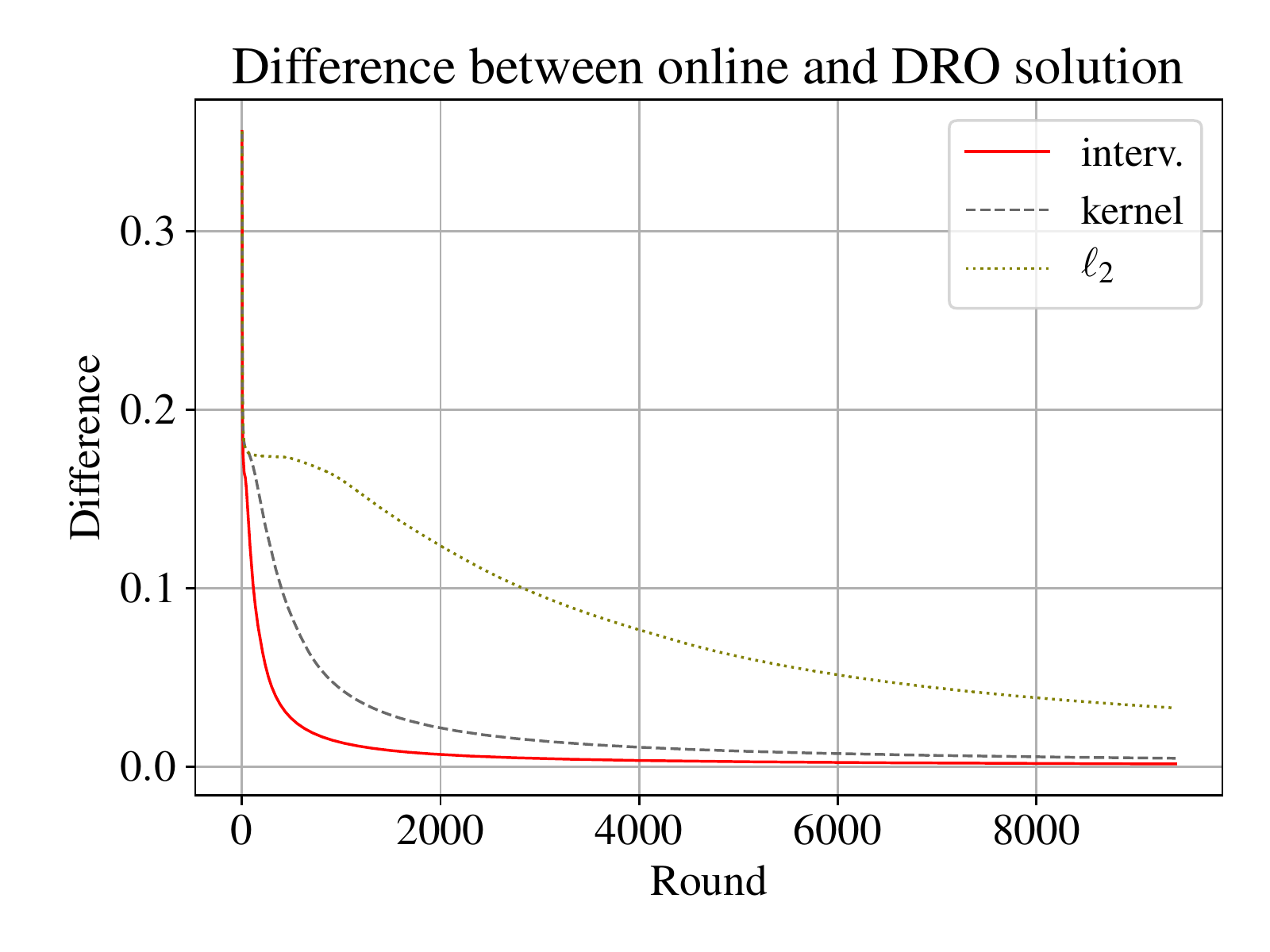}\label{fig:dynamic_regret}}
	\subfigure[Convergence of online robust and DRO solutions to the stochastic solution.]{\includegraphics[width=8cm]{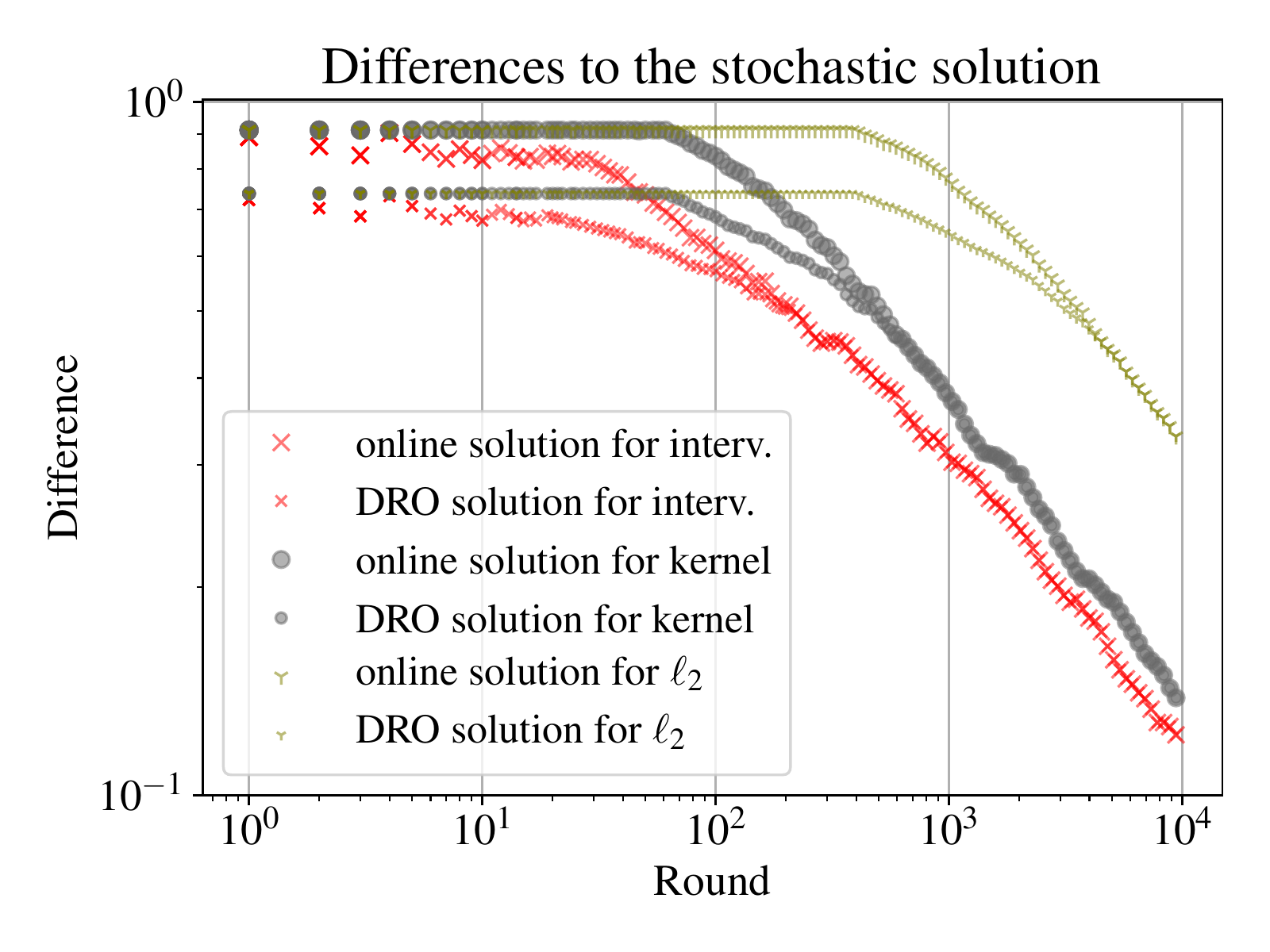}
	\label{fig:neosregret}}
	\caption{Results for \emph{blend2} with $|\mathcal{S}|=10$ and $T=10000$.}
\end{figure}

Figure~\ref{fig:mip_regret1} shows the worst-case (or guaranteed)
 expected objective value in each round for the online algorithm for $T=10000$ rounds for instance \emph{blend2} from MIPLIB with $|\mathcal{S}|=10$.
Therein, the results for all three types of ambiguity sets (confidence intervals, $\ell_2$-norm,  kernel-based) are shown under the
assumed ambiguity set in each period. 
The lines represent the corresponding average values over time. 
One can observe that the confidence intervals (red) yield the fastest
convergence and lead to solutions with least costs.
In Figure~\ref{fig:dynamic_regret} and Figure~\ref{fig:neosregret}, we plot the error between the worst-case protection of the online solutions and the DRO solutions or to the SO solutions, respectively.
It can be observed that the distance between the lines decreases over rounds. 
This means that the average error in solving the DRO problem with Algorithm~\ref{Alg:DRO_via_adversarial} shrinks and we learn the robust solution rapidly.
Both solutions also converge to the true stochastic solution.
Since, the confidence interval ambiguity sets yield a better performance and shrink faster, we focus on interval sets for the remaining experiments.

Next, we evaluate the performance of Algorithm \ref{Alg:DRO_via_adversarial}.
Table~\ref{table:avg_times_mip} shows the average running times per
iteration over all instances for different scenarios and ambiguity
sets. 
It is obvious from Table~\ref{table:avg_times_mip} that for large and difficult problems (like
e.g. nonlinear mixed-integer optimization problems), using
Algorithm~\ref{Alg:DRO_via_adversarial} allows for significant time
savings. This becomes more and more pronounced for increasing number
of iterations, as the time savings multiply by the number of rounds. 

\addtolength{\tabcolsep}{-1pt}    
\begin{table}[htb]
	\centering
	\begin{small}
		\begin{sc}
			\begin{tabular}{lrcc}
				\toprule
				& $|\mathcal{S}|$ &\begin{tabular}{c}Online\\ Robust\end{tabular} & \begin{tabular}{c}Exact\\ DRO\end{tabular}\\
				\midrule
				MIP (I) & $10$ & 52.4{s} & 115.8{s}\phantom{$^*$}\\ 
				MIP ($\ell_2$) & $10$ & 49.4{s} & 127.5{s}\phantom{$^*$}\\ 
				MIP (K) & $10$ & 56.3{s} & 129.5{s}\phantom{$^*$}\\ 
				MIP (I) & $50$ & 57.7{s} & 176.7{s}$^*$\\ 
				MIP ($\ell_2$) & $50$ & 60.4{s} & 206.1{s}$^*$\\ 
				MIP (K) & $50$ & 67.0{s} & 244.4{s}$^*$\\ 
				MIQP (I) & $2$ & 170.2{s} & 271.4{s}$^*$\\ 
				MIQP ($\ell_2$) & $2$ & 186.3{s} & 329.5{s}$^*$\\ 
				MIQP (K) & $2$ & 188.6{s} & 359.6{s}$^*$\\
				\bottomrule
			\end{tabular}
		\end{sc}
	\end{small}
	\caption{\label{table:avg_times_mip}Average running times per iteration. I:Interval, K:Kernel. $(^*)$ The DRO problems \emph{k16x240b} and \emph{10004} could not be solved within a runtime limit of one hour.}
\end{table}
%

\addtolength{\tabcolsep}{1pt}    
In more detail, in 14 out of 15 MIP instances and in 9 out of 10 MIQP ones, Algorithm~\ref{Alg:DRO_via_adversarial} was  able to run an iteration on average significantly faster than solving reformulation \eqref{Eq:DRO_ref}.
The impact of different ambiguity sets on the solution times is negligible for our approach.
However, increasing the number of scenarios amplifies the size of the
reformulated DRO problem and thus results in challenging problems for
which the online algorithm is considerably more efficient.
The main advantage of Algorithm \ref{Alg:DRO_via_adversarial} is that we avoid solving the full~\ref{Eq:DRON} problem in each round, which can be algorithmically challenging.
Instead, only the solution of an SO and a convex projection problem is
calculated in the learning decomposition approach. The computational results show that our approach is able to generate high-quality solutions within a short running time.


\subsubsection{Comparison with Other Methods}\label{Sec:comparison_with_diferent_methods}

In this section, we compare the performance of our novel algorithms against other approaches that can solve DRO problems. 
Specifically, we focus on the ambiguity sets outlined by~\citet{esfahani2018data} (Wasserstein) and~\citet{KirschnerBogunovicJegelkaKrause2020} (Distributionally Robust Bayesian Optimization or DRBO), which leverage Wasserstein and Kernel based ambiguity sets respectively.
We implement the Wasserstein and DRBO approaches as outlined in the respective papers with an exact computation of the corresponding distributionally-robust optimization problem. 
In the numerical experiments, we observe that the interval ambiguity
sets yield a solution with a comparable worst-case expected objective
as the solutions generated by the Wasserstein and DRBO methods. We
also compare our results with the worst-case performance of the stochastic optimum using the Maximum Likelihood Estimator (MLE) in each round (\emph{running SO}). 
This solution convergences to the true stochastic optimum in the limit ($T\ra \infty$) but is not protected against ambiguity and has no worst-case guarantees in contrast to distributionally robust solutions.

\begin{figure}
	\centering
	\includegraphics[width=8cm]{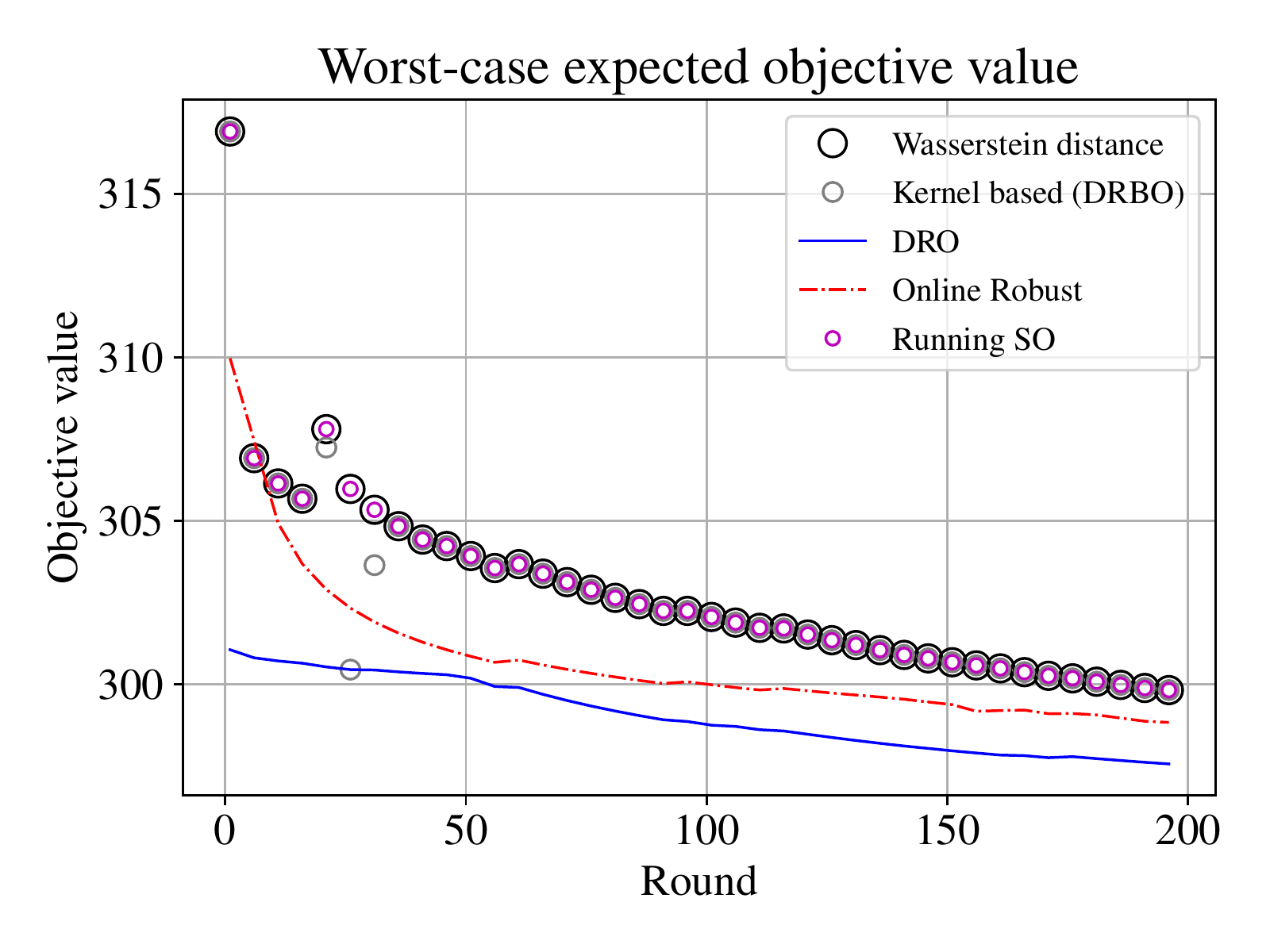}
	\caption{Results for  \emph{supportcase16} ($|\mathcal{S}|=20$, $T=100$).}
	\label{Fig:Comparison_with_other_ambiguity_sets}
\end{figure}

%
\begin{table}[h]
	\centering
	\begin{tabular}{lll}
		\toprule
		& $|\mathcal S| = 10$ & $|\mathcal S| = 50$ \\ 
		\midrule 
		DRO &	2.2s &	3.5s \\
		Wasserstein &	0.4s &	1.1s \\
		DRBO &	7.9s &	30.9s \\ \midrule
		Online robust &	0.3s & 0.4s \\
		Running SO &	0.2s &	0.2s \\
		\bottomrule
	\end{tabular}
	\caption{Avg.\ running times for  \emph{supportcase16},  $T=100$.}
	\label{Table:running_times_different_methods}
\end{table}
As a representative example, Figure~\ref{Fig:Comparison_with_other_ambiguity_sets} illustrates these results for the \emph{supportcase16} instance from MIPLIB with $|\mathcal{S}|=10$  and $T=100$. For all methods, the worst-case expected objective value shrinks over time with more data.
The online and exact solutions yield similar  worst-case protection as the other methods.
However, the online robust solution can be calculated more efficiently than the other DRO methods, cf. Table~\ref{Table:running_times_different_methods}.
Only the stochastic solutions with the MLE can be computed faster.
However they do not have solution quality guarantees under ambiguity 
and may lead to a bad worst-case objective.
Thus in total, the online robust method is preferable. 
The average
running times of all instances are given in
Table~\ref{Table:running_times_different_methods_extended} and support
the observations. Indeed, using the online robust approach, the
instances all can be solved quickly within the time limit, whereas it
takes considerably longer for the other approaches that may already
reach the time limit for some instances. 
\begin{table}[h]
	\centering
	\begin{tabular}{lccc}
		\toprule
		&\begin{tabular}{c}MIP\\ $ |\mathcal S| = 10$\end{tabular}  
		&\begin{tabular}{c}MIP\\ $ |\mathcal S| = 50$\end{tabular}
		&\begin{tabular}{c}MIQP\\ $ |\mathcal S| = 2$\end{tabular} \\ 
		\midrule 
		DRO & 45.6s\phantom{$^{**}$} & 55.9s$^*$\phantom{$^{*}$} & 271.4s$^*$\phantom{$^{*}$}  \\ 
		Wassertein & 52.3s\phantom{$^{**}$} & 59.1s\phantom{$^{**}$} & 299.9s\phantom{$^{**}$}  \\ 
		DRBO & 42.7s$^{**}$ & 66.1s$^{**}$ & 738.3s$^{*}$\phantom{$^{*}$}  \\ \midrule 
		Online robust & 26.8s\phantom{$^{**}$} & 27.1s\phantom{$^{**}$} & 170.2s\phantom{$^{**}$}  \\ 
		Running SO & 26.6s\phantom{$^{**}$} & 26.9s\phantom{$^{**}$} & 172.6s\phantom{$^{**}$} \\ 
		\bottomrule
	\end{tabular}
	\caption{Avg.\ running times for benchmark instances. ($^*$) One testcase could not be solved within one hour. ($^{**}$) Four testcases could not be solved within one hour.}
	\label{Table:running_times_different_methods_extended}
\end{table}

\subsection{Network Design under Uncertainty}
In addition to solving classical benchmark instances as performed in
the last section, we next study the solution of a practically very relevant and highly challenging
combinatorial optimization problem under uncertainy. 
Given an undirected graph $\mathcal{G}=(V,E)$, the goal of
robust network design \citep{CacchianiJuengerLiersetal.2016} is to compute a
minimal cost network topology together with the
corresponding edge capacities $f\in \Z_+^{|E|}$ in order to fulfill a given demand $b\in \Z^{|V|}$.
The demand is assumed to be uncertain and an element of the scenario set $\mathcal{S}\coloneqq \{ b_1,...,b_{|\mathcal{S}|} \} \subset \mathbb{Z}^{|V|}$ with (unknown) probability vector $p^*\in [0,1]^{|\mathcal{S}|}$.
For every edge $\{i,j\}\in E$ and every scenario $s\in \mathcal{S}$, we are given costs $c_{ijs}>0$ and flow capacities $d_{ijs},\bar{d}_{ij}>0$.
In this practical application, it can be assumed that additional
information on the demand distributions become available over time,
so that a DRO over time approach is a very natural modelling choice. 

The optimization problem for distributionally robust network design is then given by
\begin{align*}
	\min_{f}  &\text{ } \max_{p\in \mathcal{P}} \text{ } \sum_{s\in \mathcal{S}} \sum_{\{i,j\}\in E} c_{ijs} \left( f_{ijs} + f_{jis} \right)p_s \\
	\text{ s.t. } &\sum_{j:\{j,i\}\in E} f_{jis} - \sum_{j:\{i,j\}\in E} f_{ijs} = b_{is} \quad \forall i \in V, s\in \mathcal{S}, \\
	&\sum_{s\in \mathcal{S}} (f_{ijs} + f_{jis})  \leq \bar{d}_{ij} \quad \forall  \{i,j\}\in E, \\
	&f_{ijs} + f_{jis}  \leq d_{ijs} \quad \forall  \{i,j\}\in E , s\in \mathcal{S},\\
	&f_ {ijs} \in \Z_+ \quad\quad \forall \{i,j\}\in E, s\in \mathcal{S}.
\end{align*}
Here, the objective minimizes the flow costs.
The first constraint ensures that demands are satisfied. 
The second and third constraints limit the flows to be within arc capacities and the final constraint ensures integral flows. 


We evaluate the novel approach on the instances \emph{res8} ($|V|=50,
|E|=77$), \emph{w1$\_$100} ($|V|=100, |E|=207$) and \emph{w1$\_$200} ($|V|=200, |E|=775$) from
\cite{network_design_instances} and construct the scenarios as follows:
On half the nodes, we place balanced random demands from
$\{-10,...,10\}$. 
We also restrict all flow capacities to $d_{ijs}=10$,
where edge costs and coupling capacity bounds are uniformly chosen from $c_{ijs} \in
\{1,...,10\}$ and $\bar{d}_{ij}\in\{|\mathcal{S}|,...,10|\mathcal{S}|\}$. 

\begin{figure}[htb]
	\centering
	\includegraphics[width=8cm]{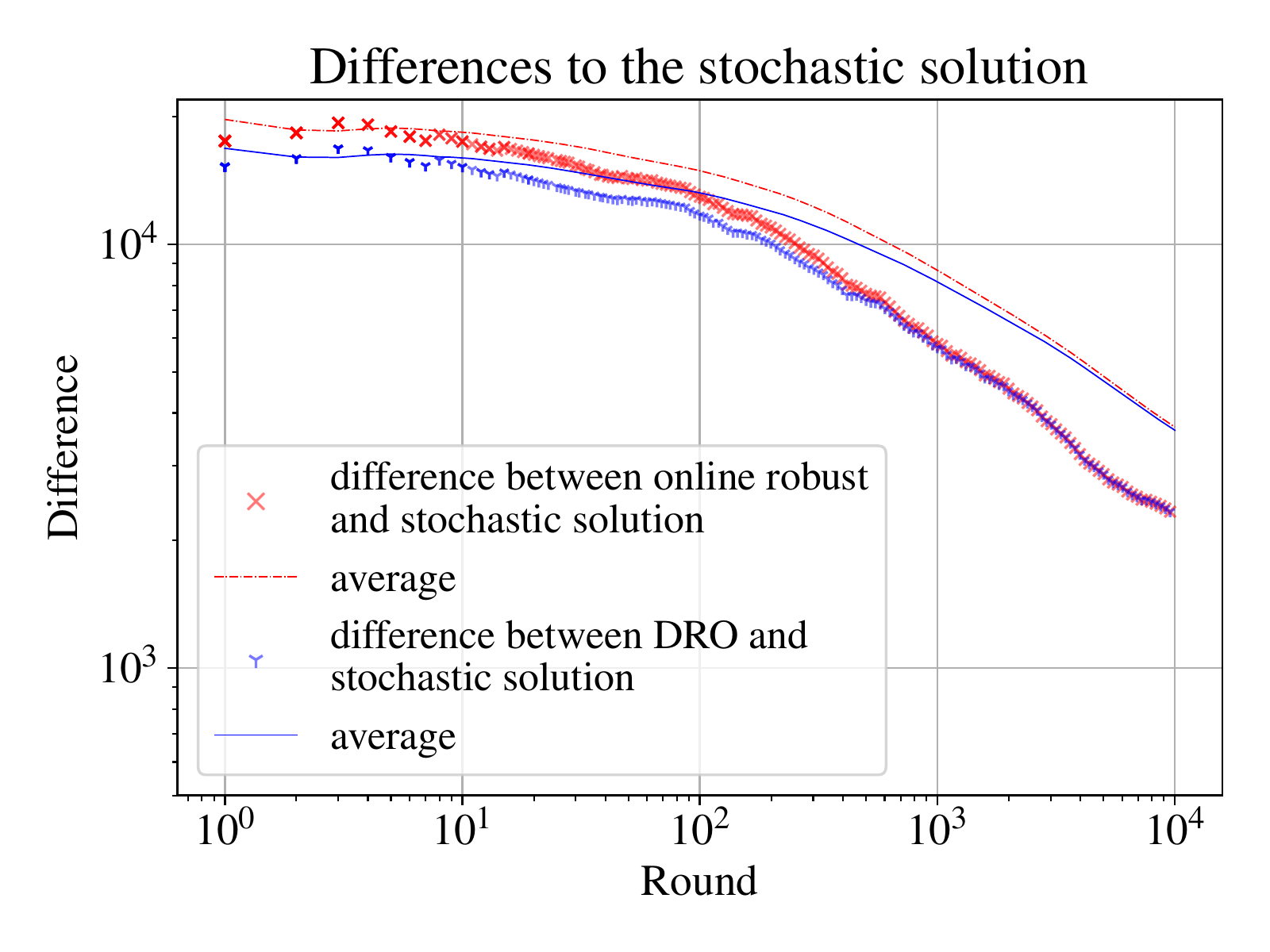}

	\caption{Results for \emph{res8} with $|\mathcal{S}|=10$ and $T = 10000$.}
	\label{fig:mip_regret1_network_design_full_regret}
\end{figure}

In Figure~\ref{fig:mip_regret1_network_design_full_regret}, the worst-case
expectation of solutions for instance \emph{res8} with $|\mathcal{S}|=10$ using confidence intervals is illustrated over $10000$ rounds on logarithmic axes. The online robust solution rapidly converges to the DRO solution.
Though starting with conservative outcomes due to limited data, the online solutions improve very quickly.
The running time benefit of the online robust approach
is clearly visible in Table~\ref{table:avg_times_network_design}.
For larger instances, it is about
272 times faster than solving reformulation \eqref{Eq:DRO_ref}.
In summary,  our method is able to generate robust solutions with shorter running time. 
\begin{table}[htb]
\centering
			\begin{sc}
				\begin{tabular}{lrcc}
					\toprule
					& $|\mathcal{S}|$&Online Robust & Exact DRO \\
					\midrule
				res8& $10$ & 0.2 {s}& 0.5 {s} \\
				res8& $50$& 0.6 {s} & 11.6 {s}\\
				w1\_100& $10$& 0.3 {s} &  32.0 {s} \\
				w1\_100& $50$& 1.5 {s} &  95.6 {s} \\
				w1\_200& $10$& 1.2 {s} &  38.7 {s} \\
				w1\_200& $50$& 4.7 {s} &  1282.2 {s} \\ \bottomrule
			\end{tabular}
			\end{sc}
	\caption{\label{table:avg_times_network_design}Average running times per iteration using confidence intervals.}
\end{table}


\subsubsection{Learning Optimal Route Choice}

In this section, we consider the problem of choosing the shortest paths in a street network
where the travel times on the arcs are affected by random deviations.
It is natural to assume that the driver gradually adapts the route
according to the observed travel times in order to reach the destination 
in the shortest possible (expected) time, making a DRO over time
approach a good modeling choice. 

This means that the driver solves the shortest-path problem
on a directed graph $ G = (V, A) $
with uncertain travel times $ c\colon A \to \R_+ $.
Let $ v_1, v_2 \in V $ be the origin and the destination, respectively.
We assume that there is a finite set
of traffic scenarios~$ {\mS}
= \{c_1, c_2, \ldots, c_{\abs{\mS}}\} \subset \R_+^{\abs{A}} $
which correspond to different realizations of the travel times on the arcs,
each materializing with an unknown probability $ p^*_k \in [0, 1] $, $ k = 1,...,\abs{\mS} $.

In each round~$ t=1, \ldots, T $,
(\eg every morning when driving to work),
the driver chooses a $ v_1 $-$ v_2 $-route
given by the vector $ x_t \in \{0, 1\}^{\abs{A}} $,
which models the edges traveled, along the path chosen, in that round.
The expected travel time in a round $ t \in T $
is then given by $ \sum_{s \in \mS} p^*_s \langle c_s, x_t \rangle $.
As the true scenario distribution is unknown, the driver is assumed to solve
the distributionally robust shortest-path problem in an online fashion,
\ie using Algorithm~\ref{Alg:DRO_via_adversarial}.

\begin{figure}
	\centering
	\subfigure[Evolution of worst-case expected objective for online robust and stochastic solutions with MLE (running SO) averaged over ten examples.]{
		\includegraphics[width=0.47\linewidth]{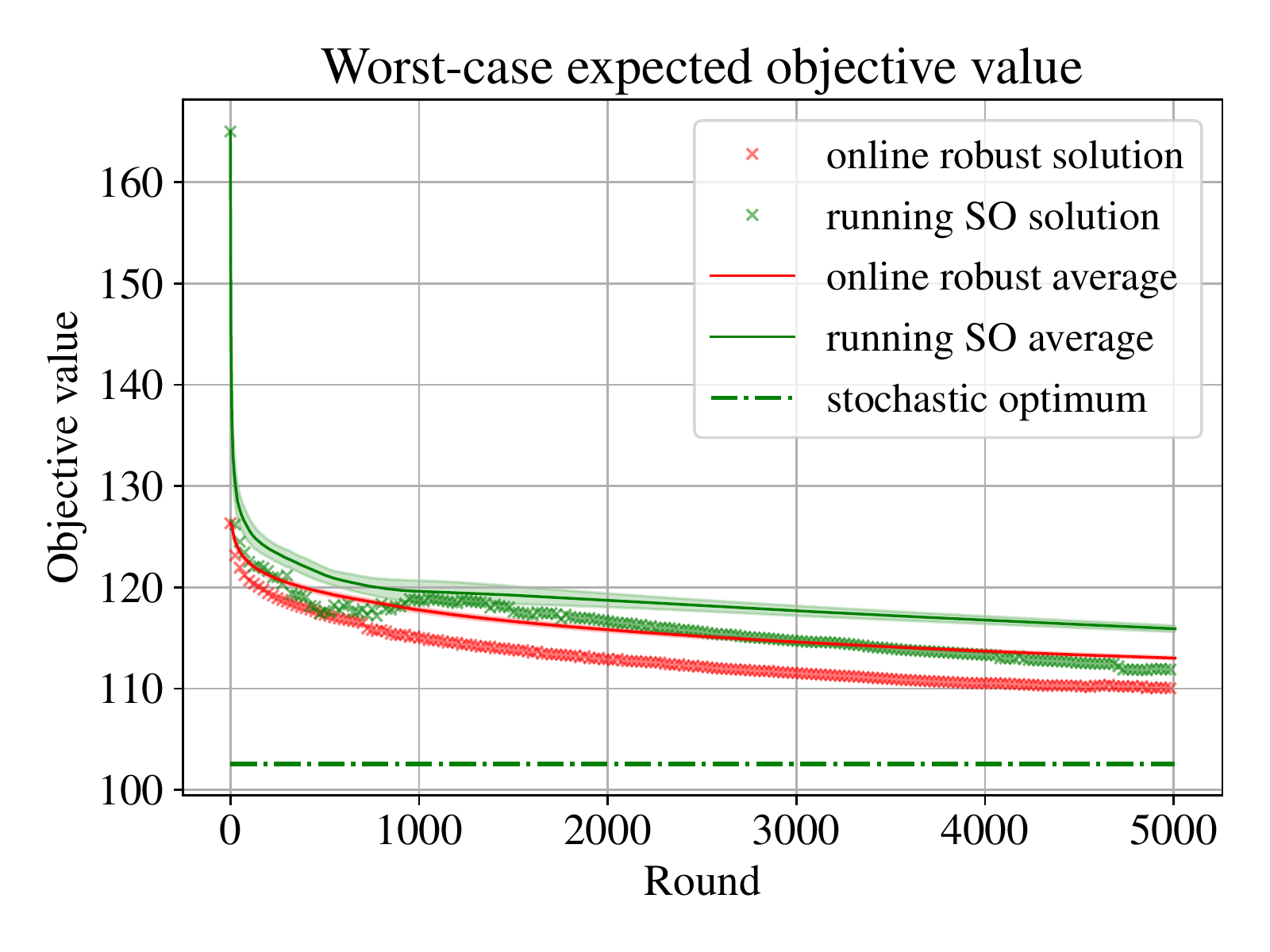}
		\label{Fig:ShortestPathProblem-22}
	}	
	\subfigure[Convergence of online robust,  DRO solutions and running SO solutions averaged over ten examples.]{
		\includegraphics[width=0.47\linewidth]{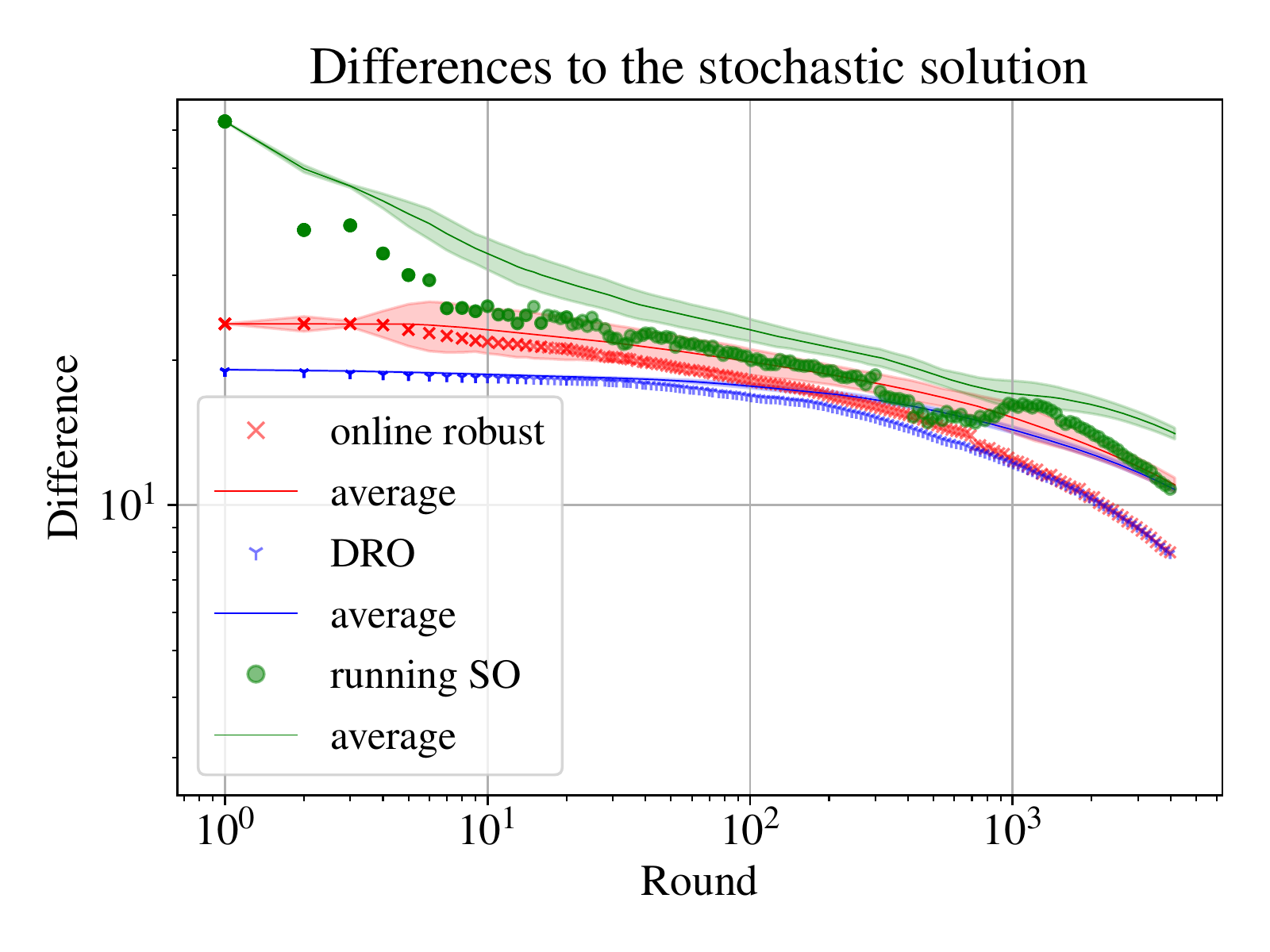}
		\label{Fig:ShortestPathProblem-21}
	}
	\caption{The outcome of Algorithm~\ref{Alg:DRO_via_adversarial} for learning an optimal route choice in terms of solution quality over time for $T=5000$.}
	\label{Fig:ShortestPathProblem}
\end{figure}
\begin{figure}
	\subfigure[Evolution of worst-case expected objective for a single example.]{
		\includegraphics[width=0.47\linewidth]{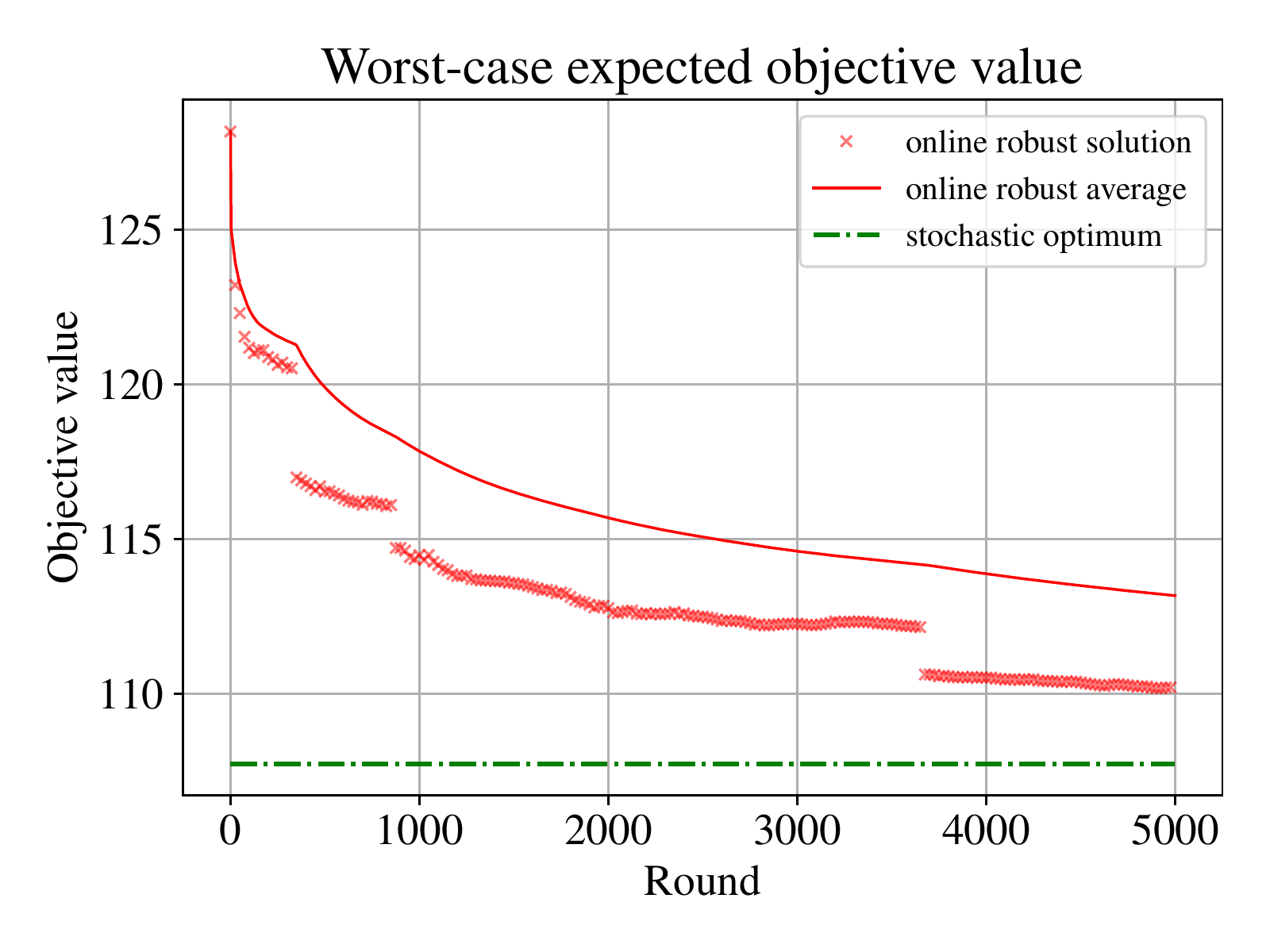}
		\label{Fig:ShortestPathProblem-SC}
	}
	\subfigure[Convergence of online robust and DRO solutions for a single example.]{
		\includegraphics[width=0.47\linewidth]{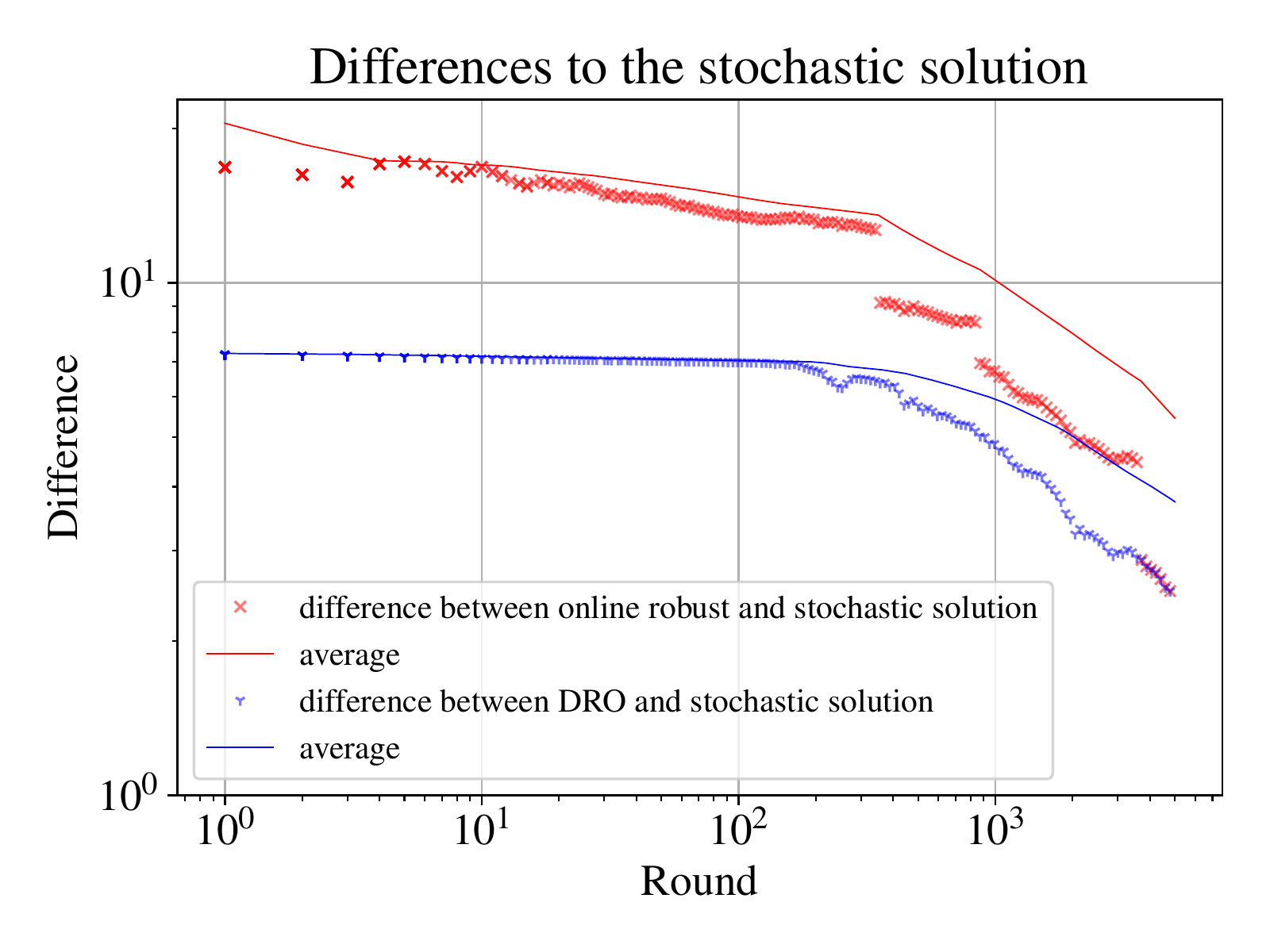}
		\label{Fig:ShortestPathProblem-Original}
	}
	
	\caption{The outcome of Algorithm~\ref{Alg:DRO_via_adversarial} for single runs that correspond to the routes given in Figure~\ref{Fig:ShortestPathProblem-Paths}.\\}
	\label{Fig:ShortestPathProblemb}
\end{figure}
In the following, we analyse the outcome of this experiment
on an aggregated version of the real-world city network of Chicago.
It is available as instance \emph{ChicacoSketch}
in Ben Stabler's library of transportation networks \citep{bstabler}
and has 933~nodes and 2950~arcs
(of which we ignore the 387~nodes representing \qm{zones}
as well as their incident arcs).
In this data set, each arc~$a$ has a certain free-flow time $ c_{\text{free}, a} $,
which we assume to be the uncongested travel time.
In addition, we generate nine congestion scenarios
by perturbing $ c_{\text{free}, a} $.
We first choose~$ v_1 $ and~$ v_2 $ such that the driven path
would span the entire extract of the city map.
Now, for all arcs~$a$ we uniformly draw $ c_{s, a} \sim [0,  2c_{\text{free}, a}]$.
Finally, we uniformly draw a random \qm{true} probability distribution~$ p^* \sim \mP$.

For the above setup, we use Algorithm~\ref{Alg:DRO_via_adversarial}
in order to let the driver iteratively adapt
to the dynamically changing travel times.
In Figure~\ref{Fig:ShortestPathProblem}, we illustrate solution quality over time.
In order to show the stability of our algorithm, we additionally repeated the experiments ten times and plot their mean solution quality as well as their standard deviation.
We observe that the online average and DRO average jointly converge
towards the expected value of the minimum expected travel time.
The regret tends to zero over the long run.

As an illustrative example, in Figures~\ref{Fig:ShortestPathProblem-SC}~and~\ref{Fig:ShortestPathProblem-Original}, one specific run is plotted to evaluate it in more detail. 
In these plots, one can see that in each of the rounds~1,~349,~860 and~3661, long stretches of the chosen path are abruptly improved. 
There are also visible jumps in the online robust solution quality.

\begin{figure}
	\phantomcaption
	\centering
	\subfigure[True stochastic solution]{
		\includegraphics[width=0.47\linewidth,trim={12cm 12cm  12cm 12cm},clip]{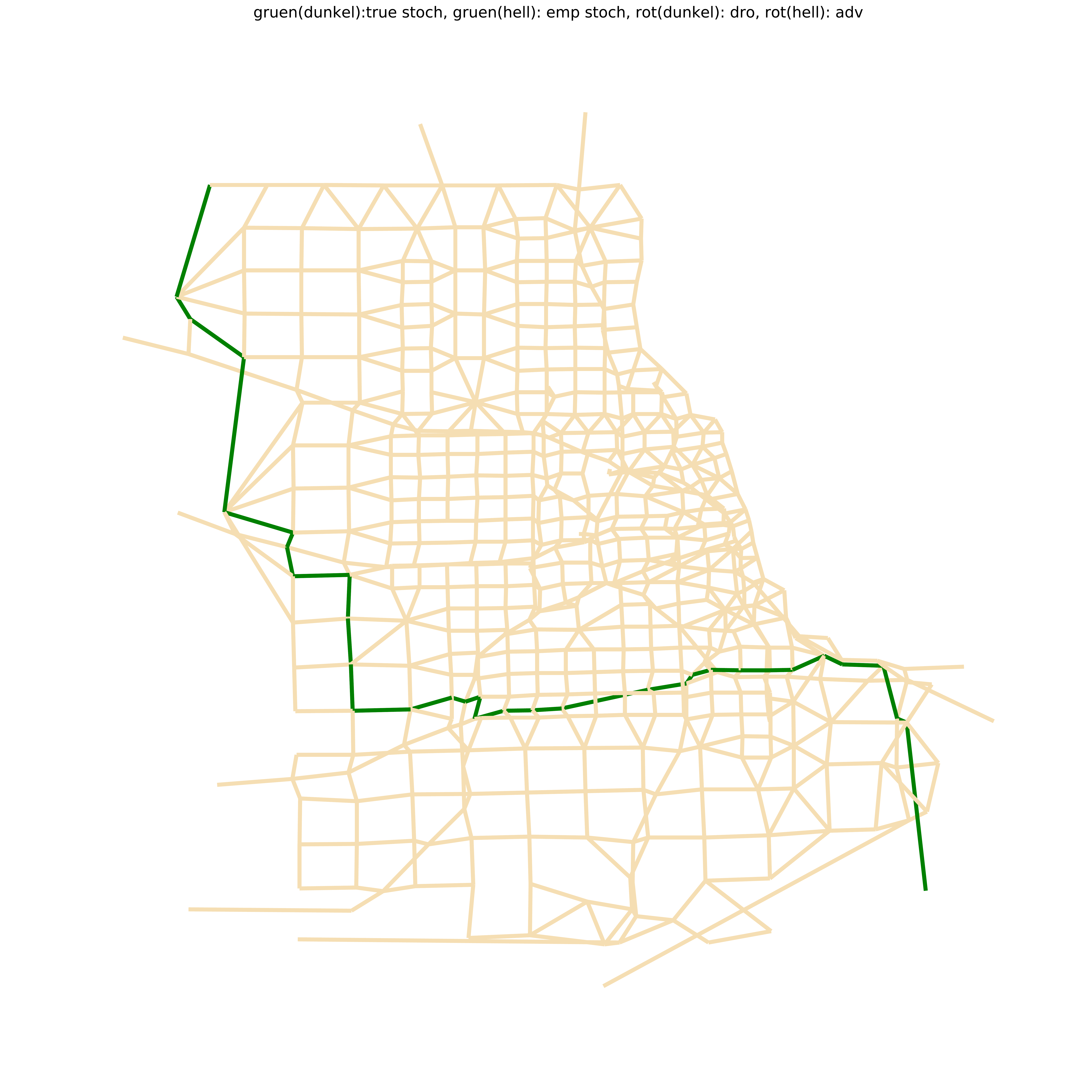}
		\label{Fig:ShortestPathProblem-Paths-Round-True0}
	}	
	\subfigure[Round 0]{
		\includegraphics[width=0.47\linewidth,trim={12cm 12cm  12cm 12cm},clip]{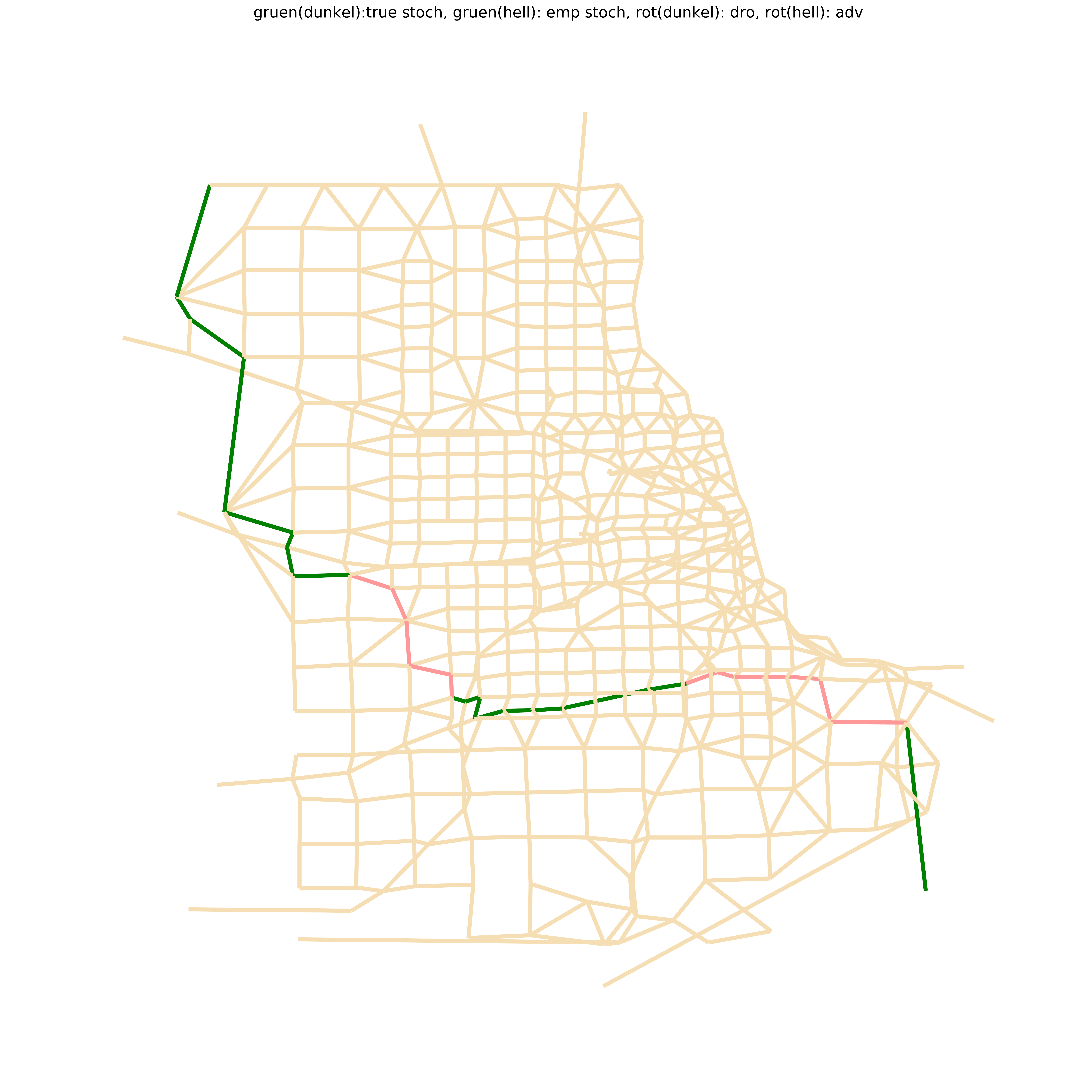}
		\label{Fig:ShortestPathProblem-Paths-Round0}
	}
	\subfigure[Round 1]{
		\includegraphics[width=0.47\linewidth,trim={12cm 12cm  12cm 12cm},clip]{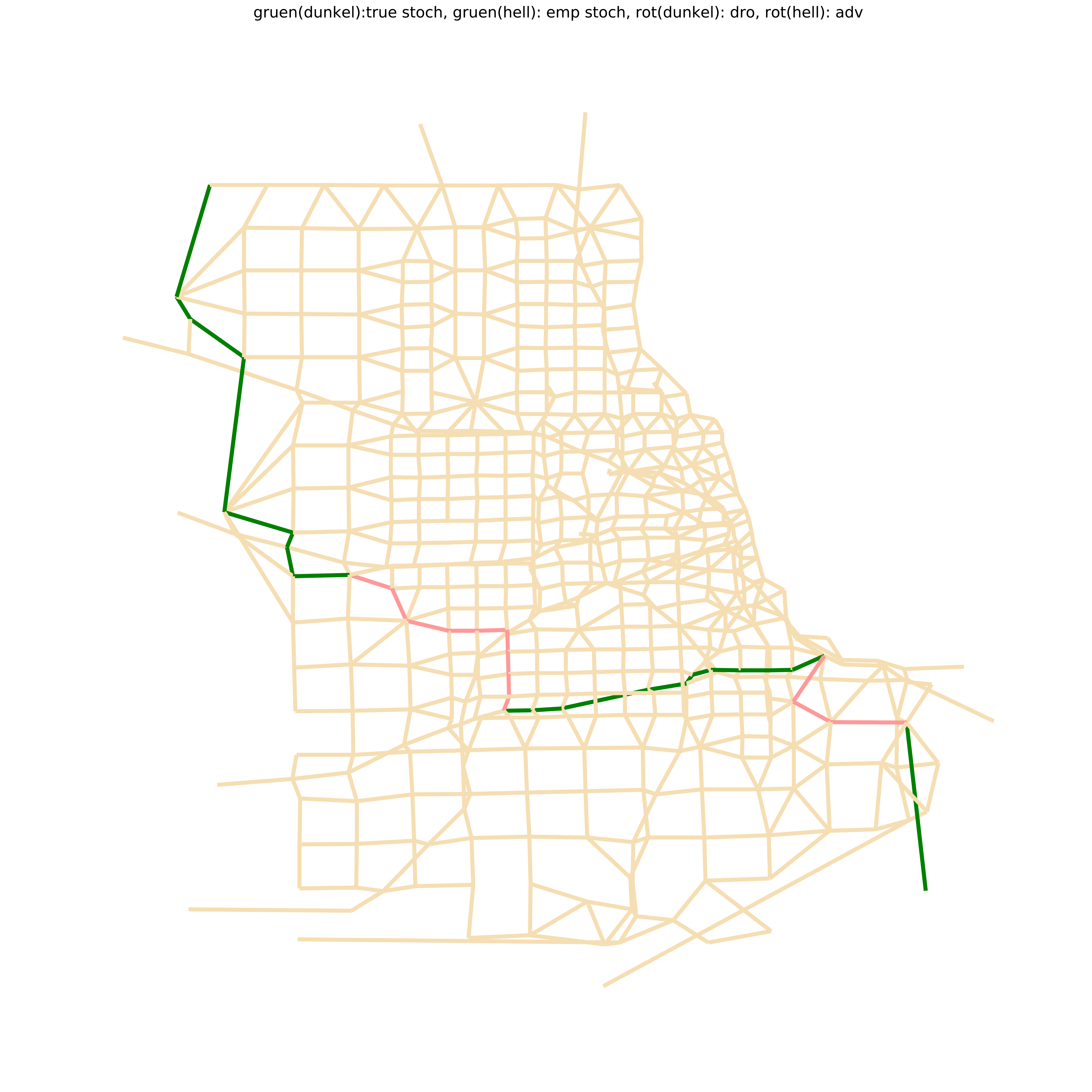}
		\label{Fig:ShortestPathProblem-Paths-Round1}
	}
\end{figure}
\begin{figure}\ContinuedFloat
	\centering
	\subfigure[Round 349]{
		\includegraphics[width=0.47\linewidth,trim={12cm 12cm  12cm 12cm},clip]{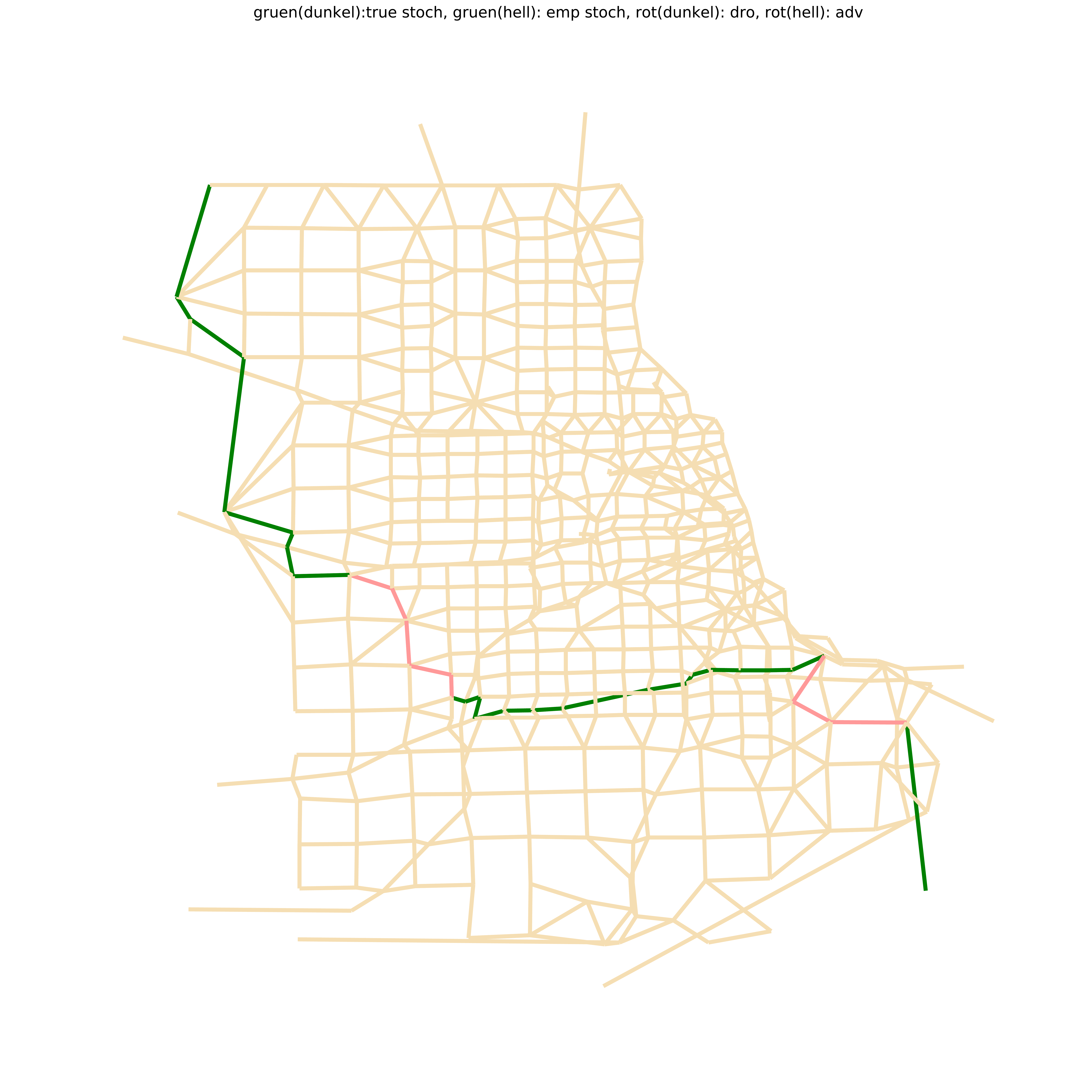}
		\label{Fig:ShortestPathProblem-Paths-Round349}
	}
	\subfigure[Round 1879]{
		\includegraphics[width=0.47\linewidth,trim={12cm 12cm  12cm 12cm},clip]{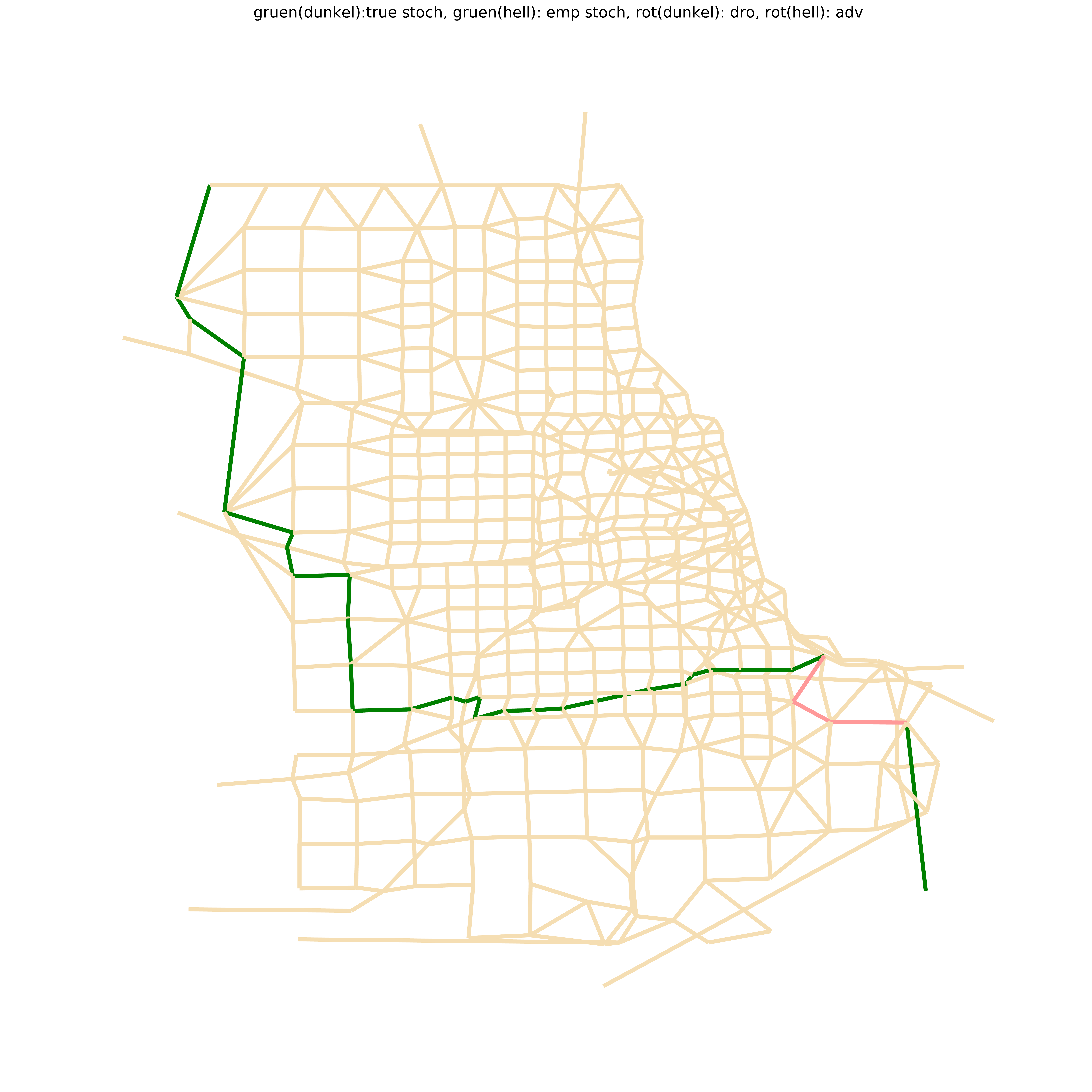}
		\label{Fig:ShortestPathProblem-Paths-Round1879}
	}
	\subfigure[Round 3661]{
		\includegraphics[width=0.47\linewidth,trim={12cm 12cm  12cm 12cm},clip]{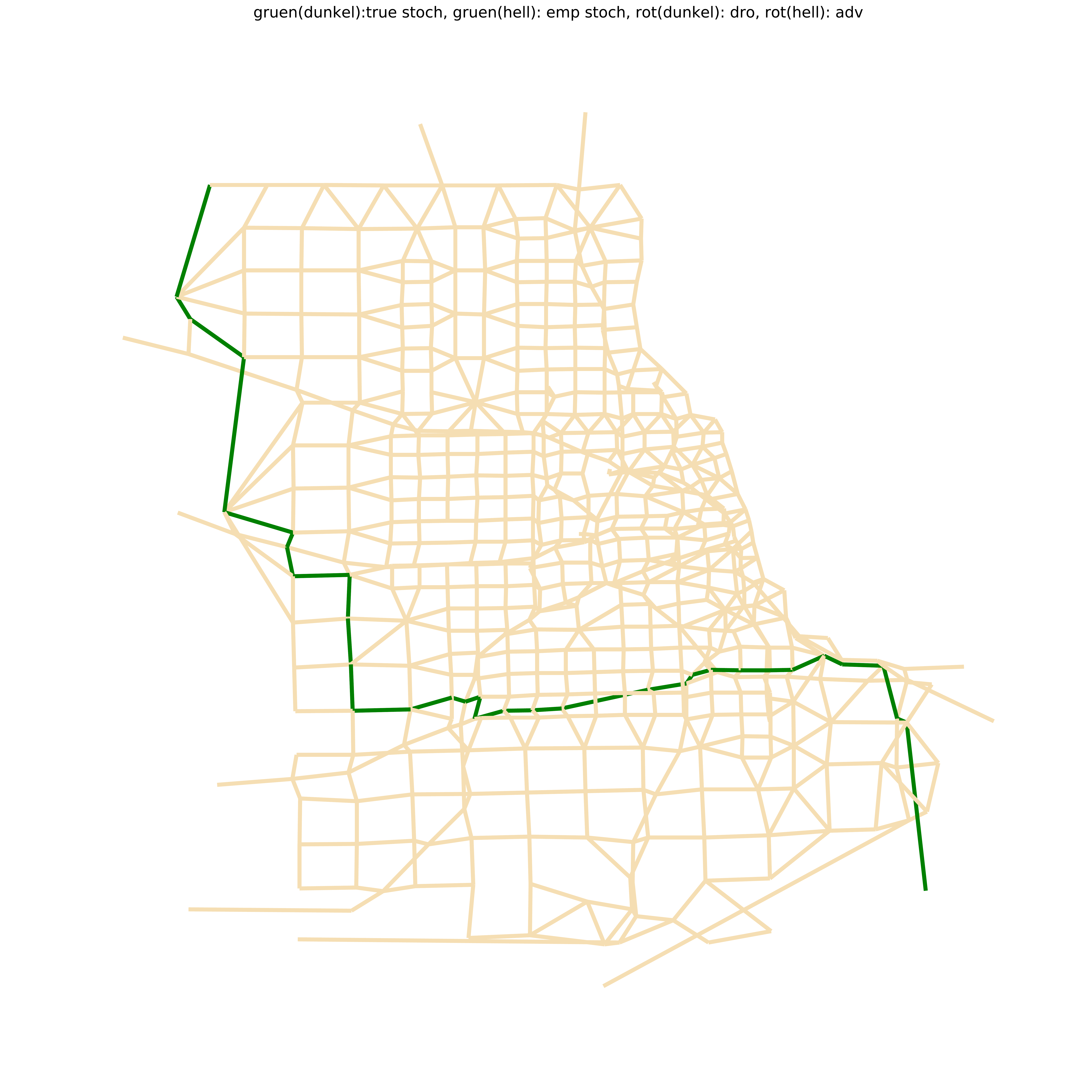}
		\label{Fig:ShortestPathProblem-Paths-Round3661}
	}
	\caption{Map (a) shows the path of the true stochastic solution in green. Maps (b)--(f) show the path taken by the online solution in those rounds in which the solution changes. If parts of the path in these pictures coincide with the path of the true stochastic solution, it is shown in green, otherwise the difference is shown in pink. The pictures show the convergence of the online solution to the true stochastic solution as time moves on.}
	\label{Fig:ShortestPathProblem-Paths}
\end{figure}

In Figure~\ref{Fig:ShortestPathProblem-Paths},
we depict how route choice evolve over the rounds $ t \in T $;
introducing a new picture whenever a structurally new solution is found.
At the beginning, the driver takes the nominally optimal path,
\ie, the optimal path w.r.t. the unperturbed cost vector.
We see two subtours in which the path differs from the true stochastic solution.
In the following iterations, there remain two deviating subtours,
but they change slightly in round~1 and round~349.
In round~1879, one of these subtours disappears,
and after round~3661 the solution coincides
with the true stochastic optimum.

Altogether, this example shows how DRO can be used to improve performance
in the face of uncertainty by leveraging information arriving over time.
At the beginning, with no or little information available,
hedging against uncertainty necessarily means to implement conservative solutions.
However, as more information on the uncertainty is gathered over time,
the solution quality improves as protection against uncertainty is
less costly. 


%% file: conclusion.tex
\section{Conclusion}\label{sec:conclusion}

We introduce a novel method for decision-making
under uncertainty
over time,
employing a combination of distributionally robust optimization
and online learning.
In each iteration, our algorithm solves a stochastic optimization problem in combination with the online gradient descent algorithm.
We show that our online algorithm converges to the exact solution of the DRO problem with an increasing amount of iterations. 
We also show that DRO solution converges to the true SO optimum in the limit.
The theoretical and numerical results
demonstrate the effectiveness of this method.
Indeed, it obtains high-quality robust solutions
with short computational times. 
Though, our work is tailored to discrete scenarios of finite dimension
and requires solving the outer minimization problem exactly, it can be applied to a wide range of practical problems with varying uncertainty models
and objective function structures.
Furthermore, our flexible framework can be extended,
either by incorporating more stochastically expressive ambiguity sets (e.g. for continuous distributions)
or using online methods also for the decision problem of the $x$-player.

%% file: extension.tex

\section*{Appendix}
Here we present the proofs omitted from the main body of the paper as well as additional numerical experiments. 

\subsection{Proofs of Theorem~\ref{thm:asymp_consis} and Theorem~\ref{thm:soln_conv}}\label{Appendix:proof_theorem_2_1_and_2_2}
In the following, we prove that the sequence of solutions to the problem~\eqref{Eq:DRON} converge to the true stochastic optimization~\eqref{Eq:Stoch_nominal} problem.
These results are an adaptation of the proofs presented in~\citet{MohajerinEsfahani2018} for our setting. 

\begin{lemma}
	\label{thm:dist_incl}
	Given the true distribution $p^*$ and the ambiguity set $\mP_t$ at any time $t$, we have
	\[\mathbb{P}\left[p^* \in \mP_t\right] \geq 1 - \delta_t \quad \text{ for all }  t=1,\dots,T.\] 
\end{lemma}
\proof{Proof:}
	This is true by the construction of the set $\mP_t$, which is such that it contains the true distribution with probability at least $1-\delta_t$. \Halmos

\begin{lemma}[Finite sample guarantee]
	Given a solution $x_t$ to the problem~\ref{Eq:DRON}, we prove that 
	\label{thm:finite_guaran}
	\[\mathbb{P}\left[\mathbb{E}_{s \sim p^*}[f(x_t, s)\right] \leq \widehat{J}_t] \geq 1 - \delta_t~\text{for all}~t=1,\dots,T.\]
\end{lemma}
\proof{Proof:}
	From Lemma~\ref{thm:dist_incl}, we know that $p^* \in \mP_t$ with probability at least $1-\delta_t$.  
	Thus, we have
	\[
	\mathbb{E}_{s \sim p^*}[f(x_t, s)] \leq \max_{p \in \mP_t} \mathbb{E}_{s \sim p} [f(x_t,s)], \]
	with probability at least $1-\delta_t$. 
	
	The right hand side (RHS) term in the above equation is the definition of $\widehat{J}_t$.
	Thus,
	\[\mathbb{E}_{s \sim p^*}[f(x_t, s)] \leq \widehat{J}_t,\]
	with probability of at least $1 - \delta_t$. \Halmos

\begin{lemma}[Borel-Cantelli Lemma]
	\label{lem:bclemma}
	Let $ E_1, E_2, \ldots $ be a sequence of events.
	If $ \sum_{i = 1}^{\infty} P(E_i) < \infty $ then 
	\[P[\text{an infinite number of } E_i \text{ occur}] = 0.\]
\end{lemma}

\begin{lemma}[Convergence of Distributions]
	\label{lem:conv_of_distr}
	Given the ambiguity set $\mP_t$, we prove that 
	\[\lim_{t \ra \infty} \sup_{p \in \mP_t} \|p - p^*\|_2  = 0 \text{ with probability } 1.\]
\end{lemma}
\proof{Proof:}
	From Lemma~\ref{Lemma:shrinking_interval},~\ref{lemma:l2_shrinking_difference} and~\ref{lemma:kernel_shrinking_difference} we know that for any of the three given types of ambiguity sets there exists a function $r(t)$ which satisfies \[\Prob \left[\sup_{p \in \mP_t}\|p - p^*\|_2 \leq r(t)\right] \geq 1- \delta_t,\]
	and $\lim_{t \rightarrow \infty}r(t) = 0$. 
	\\
	This means that 
	$$\Prob \left[\sup_{p \in \mP_t}\|p - p^*\|_2 -  r(t) > 0\right] \leq  \delta_t.$$
	By construction, it follows that $\sum_{t=1}^{\infty}\delta_t  < \infty$ (as $\delta_t = \frac{6 \delta}{\pi^2 t^2}$). 
	Then the Borel-Cantelli Lemma~\ref{lem:bclemma} implies that
	$$\Prob \left[\lim_{t \ra \infty} \sup_{p \in \mP_t}\|p - p^*\|_2 - r(t) \leq 0\right] = 1.$$
	Since $\lim_{t \ra \infty} r(t) = 0$ and $\|p - p^*\|_2  \geq 0$, this means that 
	\[\lim_{t \ra \infty} \sup_{p \in \mP_t} \|p - p^*\|_2 = 0 \text{ with probability } 1.\] \Halmos

\proof{Proof of Theorem~\ref{thm:asymp_consis}}
	We know that $x_t \in \mX$  and \\$J^* \leq \mathbb{E}_{s \sim p^*}[f(x_t, s)]$ as $x_t$ is a suboptimal solution.
	Applying Lemma~\ref{thm:finite_guaran}, we obtain
	\[
	\Prob \left[J^* \leq \mathbb{E}_{s \sim p^*}[f(x_t, s)] \leq \widehat{J}_t\right] \geq \Prob \left[p^* \in \mP_t\right] \geq 1 - \delta_t.
	\]
	Since $\sum_{t=1}^{\infty} \delta_t < \infty$, by the Borel-Cantelli lemma,  
	\[
	\Prob \left[J^* \leq \lim_{t \ra \infty}\mathbb{E}_{s \sim p^*}[f(x_t, s)\right] \leq \lim_{t \ra \infty}\widehat{J}_t] = 1.
	\]
	Let $\gamma \geq 0$.
	Since $\mX$ is compact,  there exists a $\gamma$-optimal solution $x^\gamma$ to the stochastic problem, i.e., 
	\[\mathbb{E}_{s \sim p^*} [f(x^\gamma, s)] \leq J^* + \gamma.\]
	Let $p_t^\gamma  \in \mP_t$ be a $\gamma$-optimal distribution to $x^\gamma$, i.e. 
	\[\sup_{p \in \mP_t} \mathbb{E}_{s \sim p} [f(x^\gamma, s)] \leq \mathbb{E}_{s \sim p_t^\gamma} [f(x^\gamma, s)] + \gamma.\]
	Then we can write
	\begin{align*}
	&\limsup_{t \ra \infty} \widehat{J}_t \\
	\leq &\limsup_{t \ra \infty} \sup_{p \in \mP_t} \mathbb{E}_{s \sim p}[f(x^\gamma, s)] \\
	\leq &\limsup_{t \ra \infty} \mathbb{E}_{s \sim p_t^\gamma}[f(x^\gamma, s)] + \gamma \\
	= &\limsup_{t \ra \infty} \mathbb{E}_{s \sim p^*}[f(x^\gamma, s)] + \sum_{s \in \mS} f(x^\gamma, s)(p_{st}^\gamma - p_s^*) + \gamma \\
	\leq & \limsup_{t \ra \infty} \mathbb{E}_{s \sim p^*}[f(x^\gamma, s)] + G \|p_{t}^\gamma - p^*\|_2 + \gamma \\
	= & \mathbb{E}_{s \sim p^*}[f(x^\gamma, s)] + \gamma \quad \text{ w.p. } 1\\
	= & J^* + 2\gamma \quad \text{ w.p. } 1,
	\end{align*}
	where the first inequality holds because of the definition of $\widehat{J}_t$, the second inequality holds because of the definition of $p_t^\gamma$, the first equality holds as we add and subtract $p^*$ and the third inequality holds as $|f(x,s)| \leq G \text{ for all } (x,s) \in \mX \times \mS$ by assumption.
	The final two equalities hold because of Lemma~\ref{lem:conv_of_distr} and the definition of $x^\gamma$ respectively. 
	
	With this we conclude that $\limsup_{t \ra \infty} \widehat{J}_t \leq J^*$.
	Along with the earlier assertion of $J^* \leq \lim_{t \ra \infty} \widehat{J}_t$, we can now complete the proof that $\widehat{J}_t \ra J^*$ via the sandwich argument. \Halmos

\proof{Proof of Theorem~\ref{thm:soln_conv}}
	Let $\{s_t\}_{t=1}^{\infty}$ be any sequence of scenario realizations such that 
	$\lim_{t \ra \infty}\widehat{J}_t = J^*$.
	By Theorem~\ref{thm:asymp_consis}, we have $J^* \leq \mathbb{E}_{s \sim p^*} [f(x_t, s)] \leq \widehat{J}_t$  with probability $1$. 
	By the same theorem, we also know that $\lim_{t \ra \infty} \widehat{J}_t = J^* \text{ w.p. } 1$.
	Then, we can write 
	\begin{equation}
	\label{eq:acc_obj_bound}
	\liminf_{t \ra \infty} \mathbb{E}_{s \sim p^*} [f(x_t, s)] \leq \liminf_{t \ra \infty} \widehat{J}_t = J^*.
	\end{equation}
	Consider any limit point of the sequence $\{x_t\}_{t=1}^{\infty}$. 
	Since the set $\mX$ is compact, then there exists a limit point of $\{x_t\}_{t=1}^{\infty}$ which lies in $\mX$.  
	WLOG let $x^*$ be that point and \\$\liminf_{t \ra \infty} x_t = x^*$.
	\\
	Then we have
	\begin{align*}
	J^* &\leq \mathbb{E}_{s \sim p^*} [f(x^*, s)] \\
	&= \mathbb{E}_{s \sim p^*} [\liminf_{t \ra \infty}f(
	x_t, s)] \\
	&= \sum_{s \in \mS} \liminf_{t \ra \infty}f(
	x_t, s)p_s^* \\
	&= \liminf_{t \ra \infty} \sum_{s \in \mS} f(
	x_t, s)p_s^* \leq J^*,
	\end{align*}
	where the first inequality holds because $x^* \in \mX$, the second inequality holds as $\liminf_{t \ra \infty} x_t = x^*$ and because $f(x,s)$ is continuous in $x$. 
	The second equality exploits that $\mS$ is finite and the final inequality holds because of~\eqref{eq:acc_obj_bound}.
	Thus, we have $ \mathbb{E}_{s \sim p^*} [f(x^*, s)]  = J^*$ which completes the proof. \Halmos 

\subsection{Proofs of Dynamic regret bounds}\label{App:Dynamic_regret}
In order to prove the dynamic regret bound of Theorem \ref{Th:Regret_bound}, we first show that the ambiguity sets are shrinking at a rate of $\OO(\sqrt{\log T / T})$ and contain the true data generating distribution with a high confidence. The latter is stated in the following Lemma.

A crucial point for a shrinking dynamic regret bound is that the ambiguity sets are shrinking over time.
Our ambiguity sets are constructed with increasing confidence probabilities for the multinomial distribution. 
We show that the ambiguity sets are shrinking even though the confidence $1-\delta_t$ is increasing ($\delta_t = \frac{6 \delta}{\pi^2 t^2}$).
\begin{lemma}\label{lemma:bounded_tail_gaussian}
	For the upper $(1-\frac{\delta_t}{2})$-percentile $z_{\frac{\delta_t}{2}}$ of the standard normal distribution with confidence update \\$\delta_t \coloneqq \frac{6 \delta}{\pi^2 t^2}$ and $\delta\in (0,1)$, it follows that
	\begin{align*}
	z_{\frac{\delta_t}{2}}^2 \leq 4 \log(\pi t),
	\end{align*}
	for all rounds $t=1,...,T$ .
\end{lemma}
\proof{Proof:}
	By the definition of the standard normal distribution, for the upper percentile $1-\frac{\delta_t}{2}$, we have the Gaussian tail bound
	\begin{align*}
	&1- \frac{\delta_t}{2} \leq e^{ -\frac{1}{2}z_{\frac{\delta_t}{2}}^2 }
	\\
	\implies \quad &\log(1-\frac{\delta_t}{2}) \leq - \frac{z_{\frac{\delta_t}{2}}^2}{2}
	\\
	\implies \quad &z_{\frac{\delta_t}{2}}^2 \leq - 2\log(1-\frac{\delta_t}{2}) 
	\\=&-2 \log\left( \frac{2\pi^2 t^2 - 6\delta}{2\pi^2 t^2} \right)\\
	 = &2 \log\left( \frac{\pi^2 t^2}{\pi^2 t^2 - 3\delta} \right),
	\end{align*}
	for all rounds $t=1,...,T.$
	Since $\pi^2t^2-3\delta \geq 1$ for all rounds $t=1,...,T$ and $\delta\in(0,1)$ (try $t=1$ and $\delta=1$), we are able write
	\begin{align*}
	z_{\frac{\delta_t}{2}}^2 \leq 2\log(\pi^2 t^2)\leq 4 \log(\pi t),
	\end{align*}
	for all $t=1,...,T$. \Halmos

We now prove that the ambiguity sets $\mP_t$ shrink at a rate of $\mathcal{O}(\sqrt{\log t /t})$ even as the confidence requirement $1-\delta_t$ increases. 

Here, we provide a proof for ambiguity sets as defined by confidence intervals. The proofs for the other sets are provided in the electronic companion. 

\begin{lemma}[Confidence Interval Sets]\label{lemma:shrinking_difference}
	The ambiguity sets $\mP_t$ derived from \eqref{Eq:Fitzpatrick} with confidence update $\delta_t \coloneqq \frac{6 \delta}{\pi^2 t^2}$ and $\delta\in (0,1)$ for all rounds $t=1,...,T$ fulfill 
	\begin{align*}
	\sup_{x\in \mP_{0}, y\in \mP_{1}}\Vert x - y \Vert \leq \sqrt{16|\mS| \log\pi}
	\text{ and }
	\sup_{x\in \mP_{t-1}, y\in \mP_{t}}\Vert x - y \Vert \leq \frac{ \sqrt{16|\mS| \log(\pi {(t-1)})}}{\sqrt{{t-1}}},
	\end{align*}
	for all $t=2,...,T$ with a probability of at least $1-\delta$.
\end{lemma}
\proof{Proof:}
As we know from Lemma \ref{Lemma:Union_bound} that $p^* \in \bigcap_{t=0,...,T}\mathcal{P}_t$ with a probability of at least $1-\delta$, we can compute for $t=1$:
\begin{align*}
\sup_{x\in \mP_{0}, y\in \mP_{1}} \Vert x-y \Vert &= \sup_{x\in \mP_{0}, y\in \mP_{1}} \Vert x-p^*+p^*- y \Vert \\
&\leq  \sup_{x\in \mP_{0}} \Vert x-p^*\Vert + \sup_{y\in \mP_{1}} \Vert p^*- y \Vert 
\\
&\leq \sup_{x, p\in \mP_{0}} \Vert x-p\Vert + \sup_{p,y\in \mP_{1}} \Vert p- y \Vert 
\\ &\leq \sqrt{2} + \sqrt{4|\mS|\log \pi} \leq \sqrt{16 |\mS|\log \pi},
\end{align*}
with a probability of at least $1-\delta$ because of $\sup_{x,y \in \mP_{0}} \Vert x-y \Vert = \sqrt{(1-0)^2 + (0-1)^2} = \sqrt{2}$  and Lemma~\ref{Lemma:shrinking_interval}.

Similarly for $t=2,...,T$:
\begin{align*}
\sup_{x\in \mP_{t-1}, y\in \mP_{t}} \Vert x-y \Vert &= \sup_{x\in \mP_{t-1}, y\in \mP_{t}} \Vert x-p^*+p^*- y \Vert\\ &\leq  \sup_{x\in \mP_{t-1}} \Vert x-p^*\Vert + \sup_{y\in \mP_{t}} \Vert p^*- y \Vert 
\\
&\leq \sup_{x, p\in \mP_{t-1}} \Vert x-p\Vert + \sup_{p,y\in \mP_{t}} \Vert p- y \Vert  
\\ &\leq 2\frac{\sqrt{4 |\mS| \log(\pi {(t-1)})}}{\sqrt{{t-1}}},
\end{align*}
with a probability of at least $1-\delta$. \Halmos

\begin{lemma}\label{Lemma:shrinking_interval}
	The ambiguity sets $\mP_t$ derived from \eqref{Eq:Fitzpatrick} with confidence update $\delta_t \coloneqq \frac{6 \delta}{\pi^2 t^2}$ and $\delta\in (0,1)$ for all rounds $t=1,...,T$ fulfill
	\begin{align*}
	\sup_{x, y\in \mP_t} \Vert x-y \Vert \leq \frac{ \sqrt{4|\mS|\log(\pi t)}}{\sqrt{t}}.
	\end{align*}
\end{lemma}
\proof{Proof:}
	We can compute using Lemma \ref{lemma:bounded_tail_gaussian}:
	\begin{align*}
	\sup_{x, y\in \mP_{t}} \Vert x-y \Vert^2 \leq \sum_{k = 1}^{|\mS|} (u_{kt} - l_{kt})^2 &= \sum_{k = 1}^{|\mS|} \frac{z_{\frac{\delta_t}{2}}^2}{t} 
	\leq \frac{4|\mS|\log(\pi t)}{t}.
	\end{align*} \Halmos

Finally, In the following result, we prove that the distance between elements from consecutive $\ell_2$-norm and kernel based ambiguity sets $\mP_{t-1}$ and $\mP_{t}$ also shrinks as measured in the $\ell_2$-norm while providing tighter guarantees of inclusion on the true distribution..

\begin{lemma}[$\ell_2$-norm Sets]{\label{lemma:l2_shrinking_difference}}
	Given an ambiguity set of the from $\mP_t = \left\{p \in \mP \mid \|p - \hat{p}\|_2 \leq \epsilon_t \right\}$ with $\epsilon_t\coloneqq\sqrt{\frac{2|\mS| \log 2/\delta_t}{t}} $ and $\delta_t = \frac{6 \delta}{\pi^2 t^2}$, {then for $|\mS| \geq 2$} we have
	$$\sup_{x \in \mP_{0}, y \in \mP_1} \|x - y\|_2 \leq 4\sqrt{|\mS| \log (\pi / \sqrt{3 \delta})}
	\text{ and }
	\sup_{x \in \mP_{t-1}, y \in \mP_t} \|x - y\|_2 \leq 4\sqrt{\frac{|\mS| \log (\pi (t-1) / \sqrt{3 \delta})}{{t-1}}},$$
	with probability at least $ 1-\delta$.
\end{lemma}
\proof{Proof:}
{
	For the case $t=1$ we have
	\begin{align*}
		\sup_{x \in \mP_{0}, y \in \mP_1} \|x - y\|_2 &\leq \sup_{x \in \mP_{0}} \|x - p^*\|_2 + \sup_{y \in \mP_1}  \|p^* - y\|_2\\
		&\leq 2 + 2\sqrt{|\mS| \log (\pi / \sqrt{3 \delta})}\\
		&\leq 4\sqrt{|\mS| \log (\pi / \sqrt{3 \delta})}. 
	\end{align*}
	The last inequality occurs because $\sqrt{|\mS| \log (\pi / \sqrt{3 \delta})} > 1$ for $|\mS| \geq 2$.
}
	Now for the case $t > 1$, given the true distribution $p^*$, we can write,
	\begin{align*}
		\|x - y\|_2 &\leq \|x - p^*\|_2 + \|p^* - y\|_2\\
		&\leq 2 \epsilon_{{t-1}}.
	\end{align*}
	{
	Thus, 
	$$
	\sup_{x \in \mP_{t-1}, y \in \mP_t} \|x - y\|_2 \leq 4\sqrt{\frac{|\mS| \log (\pi {(t-1)} / \sqrt{3 \delta})}{{t-1}}}.
	$$
	}
	Here, the first inequality arises from triangle inequality. 
	The second inequality is due to the fact that the true distribution is contained inside all sets $\mP_t$ for $t=1,\dots,T$ with probability at least $1-\delta$.  
\Halmos

\begin{lemma}[Kernel Based Sets]{\label{lemma:kernel_shrinking_difference}}
	Given an ambiguity set of the form $\mP_t = \left\{p \in \mP \mid \|p - \hat{p}\|_M \leq \epsilon_t \right\}$ with $\epsilon_t\coloneqq \frac{{\sqrt{C}}}{\sqrt{t}}(2 + \sqrt{2 \log(1/\delta_t)}) $ with $\delta_t = \frac{6 \delta}{\pi^2 t^2}$ 
	we have for $t \geq 2$, 
	\begin{align*}
		&\sup_{x \in \mP_{t-1}, y \in \mP_t} \|x - y\|_2 \leq \frac{8\sqrt{C}}{\lambda \sqrt{{t-1}}}(\sqrt{ \log(\pi {t}/\sqrt{6\delta})}),\\
		&\text{ and } \sup_{x \in \mP_0, y \in \mP_1} \|x - y\|_2 \leq 2+\frac{4\sqrt{C}}{\lambda}\;\; \text{ for } t = 1,\\
	\end{align*}
	with probability at least $1-\delta$.
\end{lemma}
\proof{Proof:}
	For $t=1$, we have
	\begin{align*}
		\|x - y\|_2 &\leq \|x - p^*\|_2 + \|p^* - y\|_2\\
		&\leq 2 + \frac{1}{\lambda} \|p^* - y\|_M\\
		&\leq 2 + \frac{4\sqrt{C}}{\lambda}.
	\end{align*}
	The second inequality comes from the definition of the set $\mP_0$ and the definition of the norm $\|\cdot\|_M$.

	For the case $t \geq 2$, we get.
	\begin{align*}
		\|x - y\|_2 &\leq \|x - p^*\|_2 + \|p^* - y\|_2\\
		&\leq \frac{1}{\lambda} (\|x - p^*\|_M + \|p^* - y\|_M)\\
		&\leq \frac{4}{\lambda}\frac{\sqrt{C}}{\sqrt{{t-1}}}(1 + \sqrt{ \log(\pi ({t-1})/\sqrt{6\delta})})\\
		&\leq \frac{4}{\lambda}\frac{\sqrt{C}}{\sqrt{{t-1}}}(1 + \sqrt{ \log(\pi {t}/\sqrt{6\delta})}).
	\end{align*}
	Here, the first inequality comes from the triangle inequality, the second from the fact that $\sqrt{x^{\top}Mx} \geq \lambda\|x\|_2$ for any positive definite matrix $M$ with the minimum eigen value $\lambda$. 
	The final inequality arises from the construction of the ambiguity set which contains the true distribution with probability at least $1-\delta$. 
	Now note that $\pi^2 t^2 / 6 \delta > 3 $ for $t \geq 2$.
	Thus, we can write
	\begin{align*}
		\|x - y\|_2 \leq \frac{8}{\lambda}\frac{\sqrt{C}}{\sqrt{{t-1}}}\sqrt{ \log(\pi {t}/\sqrt{6\delta})}.
	\end{align*}
  \Halmos

The properties of the confidence interval ambiguity sets mentioned previously and those of the $\ell_2$-norm and kernel based sets prove in Lemmas~\ref{lemma:l2_shrinking_difference} and~\ref{lemma:kernel_shrinking_difference} along with the cumulative path length bounds from Theorem~\ref{thm:path_length_bounds} enable us to prove the following shrinking dynamic regret.
%


\begin{theorem}[Dynamic regret bound]
	\label{Th:Regret_bound3}
	Let $f:\mathcal{X}\times\mathcal{S}\rightarrow\mathbb{R}$ be uniformly bounded, i.e., for all $(x,s) \in \mathcal{X}\times\mathcal{S}$, a constant $G>0$ exists such that $|f(x,s)|\leq G$.  Let $\eta :=  \sqrt{\frac{3 + 2 h'(T)}{TG^2|\mS|}}$ with $ \sum_{t=2}^{T}\|p - q\| \leq  h'(T) $ for $p \in \mP_{t-1}$ and $q \in \mP_t$ \\
	The output $(x_1,...,x_{T})$ from Algorithm with confidence update $\delta_t \coloneqq \frac{6 \delta}{\pi^2 t^2}$ and $\delta\in (0,1)$ fulfills
	\begin{alignat*}2
		\frac{1}{T} &\sum_{t=1}^{T} \left( \max_{p \in \mP_t}\mathbb{E}_{s \thicksim p}\left[f(x_t,s)\right]  - \min_{x\in \mathcal{X}} \max_{p\in \mathcal{P}_{t}} \mathbb{E}_{s \thicksim p}\left[f(x,s)\right] \right)
		 \leq G \sqrt{\frac{3|\mS| + 2 |\mS| h'(T)}{T}} + \frac{2G}{T},
	\end{alignat*}
	with probability at least $1-\delta$.
\end{theorem}

\proof{Proof of Theorem~\ref{Th:Regret_bound3}}
	Define $g_t(p)\coloneqq -\mathbb{E}_{s \thicksim p}\left[f(x_t,s)\right]$. A gradient descent iteration is given by
	\begin{align*}
		p_{t+1} = \arg\min_{p\in \mathcal{P}_t} \text{ }\left\langle \eta \nabla g_t(p_t), p \right\rangle + \frac{1}{2}\Vert p - p_t \Vert^2,
	\end{align*}
	with optimality criteria $
	\left\langle \eta \nabla g_t(p_t), u_t - p_{t+1}  \right\rangle + \left\langle p_{t+1} - p_t, u_t - p_{t+1}  \right\rangle \geq 0 $
	for all  $u_t \in \mathcal{P}_t$. Classical theory for gradient descent yields
	\begin{align*}
		\left\langle \eta \nabla g_t(p_t), p_t - u_t  \right\rangle \leq &\frac{1}{2}\Vert p_t - u_t \Vert^2 - \frac{1}{2}\Vert p_{t+1}-u_t  \Vert^2 + \frac{\eta^2}{2} \Vert \nabla g_t(p_t) \Vert ^2.
	\end{align*}
	Summation over rounds $t=1,...,T$ results in the following inequality	for all $u_t \in \mathcal{P}_t$:
	\begin{align}
		&\sum_{t=1}^T\left\langle \eta \nabla g_t(p_t), p_t - u_t  \right\rangle \leq  \sum_{t=1}^T\frac{\eta^2}{2} \Vert \nabla g_t(p_t) \Vert ^2. \label{Eq:sum_over_rounds} 
		+ 	\sum_{t=1}^T \left(\frac{1}{2}\Vert p_t - u_t \Vert^2 - \frac{1}{2}\Vert p_{t+1} - u_t \Vert^2 \right). 
	\end{align}
	Next, we rearrange the terms on the RHS as
	\begin{align*}
		&\sum_{t=1}^T \left(\frac{1}{2}\Vert p_t - u_t \Vert^2 - \frac{1}{2}\Vert p_{t+1} - u_t \Vert^2 \right)\\
		&= \frac{1}{2}\sum_{t=1}^{T}\left(\|p_t\|^2 - \|p_{t+1}\|^2\right) + \frac{1}{2}\sum_{t=1}^{T}2\left(p_{t+1}  - p_t\right)\cdot u_t \\
		&= \frac{1}{2}\sum_{t=1}^{T}\left(\|p_t\|^2 - \|p_{t+1}\|^2\right) + \sum_{t=1}^{T}\left(p_{t+1}  - p_t\right)\cdot u_t.
	\end{align*}
	Now consider the expression $\sum_{t=1}^{T}\left(p_{t+1}  - p_t\right)\cdot u_t $. We can write it as 
	\begin{align*}
		&\sum_{t=1}^{T}\left(p_{t+1}  - p_t\right)\cdot u_t \\
		&= (p_2 - p_1)\cdot u_1 + (p_3 - p_2)\cdot u_2 + \dots + (p_{T+1} - p_T)\cdot u_T\\
		%
		%
		&=  p_{T+1} \cdot u_T- p_1 \cdot u_1 + \sum_{t=2}^{T}(u_{t-1} - u_t)\cdot p_t.
	\end{align*}
	For the other expression we have 
	\begin{align*}
		\frac{1}{2}\sum_{t=1}^{T}\left(\|p_t\|^2 - \|p_{t+1}\|^2\right)  = \frac{1}{2}\|p_1\|^2 - \frac{1}{2}\|p_{T+1}\|^2.
	\end{align*}
	This then allows us to write 
	\begin{align*}
		&\sum_{t=1}^T\left\langle \eta \nabla g_t(p_t), p_t - u_t  \right\rangle \leq \sum_{t=1}^T\frac{\eta^2}{2} \Vert \nabla g_t(p_t) \Vert ^2 + \frac{1}{2}\|p_1\|^2  \\
		& - \frac{1}{2}\|p_{T+1}\|^2+ p_{T+1} \cdot u_T- p_1 \cdot u_1 + \sum_{t=2}^{T}(u_{t-1} - u_t)\cdot p_t . 
	\end{align*}
	We know that $\|p\| \leq 1$ and for any 2 probability vectors $p$ and $q$ we have that $0 \leq p\cdot q \leq 1$. Thus we can write		
	\begin{align*}
		&\sum_{t=1}^T\left\langle \eta \nabla g_t(p_t), p_t - u_t  \right\rangle \leq 
		\sum_{t=1}^T\frac{\eta^2}{2} \Vert \nabla g_t(p_t) \Vert ^2  + 	 \frac{1}{2}+ 1 + \sum_{t=2}^{T}(u_{t-1} - u_t)\cdot p_t. 
	\end{align*}
	Note that 
	\begin{align*}
		\Vert \nabla g_t(p_t) \Vert^2 = \sum_{k = 1}^{|\mathcal{S}|} |f(x_t,s_k)|^2 \leq |\mS|G^2,
	\end{align*}
	This, along with the fact that $\sum_{t=2}^{T}(u_{t-1} - u_t)\cdot p_t  \leq \sum_{t=2}^{T}\|u_{t-1} - u_t\| \| p_t\| \leq \sum_{t=2}^{T}\|u_{t-1} - u_t\|$, allows us to write
	\begin{align*}
		\sum_{t=1}^T\left\langle \nabla g_t(p_t), p_t - u_t  \right\rangle \leq &\frac{\eta T}{2}  |\mS|G^2  + \frac{3}{2\eta} 
		+ \frac{1}{\eta}\sum_{t=2}^{T}\|u_{t-1} - u_t\|. 
	\end{align*}    
	By assumption, $ \sum_{t=2}^{T}\|u_{t-1} - u_t\|  \leq  h'(T)$ with probability at least $1-\delta$ for some function $h'(\cdot)$. Then we have 
	\begin{align*}
		\sum_{t=1}^T\left\langle \nabla g_t(p_t), p_t - u_t  \right\rangle \leq \frac{\eta T}{2}  |\mS|G^2  + \frac{3}{2\eta} + \frac{h'(T)}{\eta},
	\end{align*}    
	with probability at least $1-\delta$. 
	Choosing the optimal $\eta = \sqrt{\frac{3 + 2 h'(T)}{TG^2|\mS|}}$, we can write our result as 
	\begin{align*}
		\sum_{t=1}^T\left\langle \nabla g_t(p_t), p_t - u_t  \right\rangle \leq 2 \sqrt{\frac{T}{2}  |\mS|G^2  \cdot  (\frac{3}{2} + h'(T))}.
	\end{align*}    
%
	Then we can write
	\begin{align*}
		\sum_{t=1}^T\left\langle \nabla g_t(p_t), p_t - u_t  \right\rangle &\leq G |\mS|^{\frac{1}{2}} \sqrt{3T + 2T h'(T)}. 
		%
	\end{align*}    
	Since $g_t(p)=-\mathbb{E}_{s\thicksim p}\left[f(x_t,s)\right]$ is linear in $p$ for all $t=1,...,T$, it follows
	\begin{align*}
		&\sum_{t=1}^T \left( \mathbb{E}_{s\thicksim u_t}\left[f(x_t,s)\right] - \mathbb{E}_{s\thicksim p_t} \left[f(x_t, s)\right] \right) \\
		&= \sum_{t=1}^T \left(g_t(p_t) - g_t(u_t) \right)= \sum_{t=1}^T \left\langle \nabla g_t(p_t), p_t - u_t \right\rangle 
		\\ 
		&\leq G |\mS|^{\frac{1}{2}} \sqrt{3T + 2T  h'(T)}. 
	\end{align*}
	Now we choose in each round $t=1,...,T$ the worst-case $ 
	u_t \coloneqq \arg \max_{p \in \mP_t} \mathbb{E}_{s\thicksim p}\left[f(x_t,s)\right] \in \mP_t $
	and recall $x_t = \arg\min_{x\in \mathcal{X}} \mathbb{E}_{s \sim p_t}\left[f(x,s)\right]$ to obtain 
	\begin{align*}
		\sum_{t=1}^T &\left( \mathbb{E}_{s\thicksim u_t}\left[f(x_t,s)\right] - \mathbb{E}_{s\thicksim p_t} \left[f(x_t, s)\right] \right) 
		= \sum_{t=1}^T \left( \max_{p \in \mP_t}\mathbb{E}_{s\thicksim p}\left[f(x_t,s)\right] - \min_{x \in \mX} \mathbb{E}_{s\thicksim p_t} \left[f(x, s)\right] \right).
	\end{align*}
	Since $p_t \in \mP_{t-1}$, we know that $\min_{x \in \mX} \mathbb{E}_{s\thicksim p_t} \left[f(x, s)\right] \leq \min_{x \in \mX} \max_{p \in \mP_{t-1}}\mathbb{E}_{s\thicksim p} \left[f(x, s)\right]$ for all $t=1,...,T$  and thus we can conclude
	\begin{align*}
		\sum_{t=1}^T& \left( \max_{p \in \mP_t} \mathbb{E}_{s\thicksim p}\left[f(x_t,s)\right] - \min_{x \in \mX} \max_{p \in \mP_{t-1}} \mathbb{E}_{s\thicksim p} \left[f(x, s)\right] \right)
		\leq G |\mS|^{\frac{1}{2}} \sqrt{3T + 2 Th'(T)}, 
	\end{align*}
	with a probability of at least $1-\delta$.
	We add and subtract $\min_{x \in \mX} \max_{p \in \mP_t} \mathbb{E}_{s\thicksim p} \left[f(x, s)\right]$ on the LHS.  Rearranging the terms this allows us to write the LHS as
	\begin{align*}
		&\sum_{t=1}^T \left( \max_{p \in \mP_t} \mathbb{E}_{s\thicksim p}\left[f(x_t,s)\right] - \min_{x \in \mX} \max_{p \in \mP_t} \mathbb{E}_{s\thicksim p} \left[f(x, s)\right] \right)\\ & + \sum_{t=1}^T \min_{x \in \mX} \max_{p \in \mP_t} \mathbb{E}_{s\thicksim p} \left[f(x, s)\right]  -  \min_{x \in \mX} \max_{p \in \mP_{t-1}} \mathbb{E}_{s\thicksim p} \left[f(x, s)\right]. 
	\end{align*}
	Observing that the last two terms telescope, bringing them to the RHS and using the upper bound $G$ on $|f(x,s)|$ we can conclude
	\begin{align*}
		\sum_{t=1}^T& \left( \max_{p \in \mP_t} \mathbb{E}_{s\thicksim p}\left[f(x_t,s)\right] - \min_{x \in \mX} \max_{p \in \mP_t} \mathbb{E}_{s\thicksim p} \left[f(x, s)\right] \right)
		\leq G |\mS|^{\frac{1}{2}} \sqrt{3T + 2 Th'(T)}  + 2G \;\;\;\;\;\text{ w.p. } 1 - \delta.
	\end{align*}
	Dividing by $T$ on both sides completes the proof. \Halmos

\begin{theorem}
	\label{thm:path_length_bounds}
	Given ambiguity sets of the form specified in Section~\ref{sec:dddro}, we have 
	\begin{align*}
	\frac{1}{2}\sum_{t=1}^{T}\Vert p_t - q_t \Vert^2 \leq h(T) &\text{ and }   \sum_{t=2}^{T}\|p_t - q_t\| \leq  h'(T)
	\text{ for all } p_t \in \mP_{t-1}, q_t \in \mP_t,
	\end{align*}
	with probability at least $1-\delta$.
	The functions $h(T)$ and $h'(T)$ for different categories of ambiguity sets are as given below: 
	\begin{enumerate}
		\item \textit{Confidence Intervals}: 
		$$h(T) = 8 |\mS| \log(\pi T) ({2} + \log T)$$
		$$h'(T) = 8 \sqrt{|\mS|T\log(\pi T)}$$
		\item \textit{Kernel based ambiguity sets, where $\lambda$ denotes the smallest eigenvalue of the kernel matrix $M$:}
		$$h(T) =  {\frac{1}{2}\left(2+ \frac{4\sqrt{C}}{\lambda}\right)^2}+ \frac{32C}{\lambda^2 } \log(\pi T/\sqrt{6\delta})(1 + \log T)$$
		\begin{align*}
			h'(T) = &\frac{16\sqrt{C}}{\lambda} \sqrt{T \log(\pi T/\sqrt{6\delta})} 	
		\end{align*}
		\item \textit{$\ell_2$-norm ambiguity sets:}
		$$
		h(T) = 8 |\mS|\log \frac{\pi T}{\sqrt{3 \delta}} ({2} + \log T)
		$$
		$$
		h'(T) = 8\sqrt{|\mS|T\log\frac{\pi T}{\sqrt{3\delta}}}
		$$
	\end{enumerate}
\end{theorem}
\proof{}
	\textbf{Confidence Intervals}.
	Calculating the function $h(T)$, we have 
	\begin{align*}
		\frac{1}{2} \sum_{t=1}^{T}\|p_t - q_t\|^2 &\leq  \frac{1}{2} 16 |\mS| \log\pi + \frac{1}{2} \sum_{t=2}^{T}16 |\mS| \frac{\log(\pi (t-1))}{t-1}\\
		&\leq 8 |\mS| \log\pi + 8 |\mS| \log(\pi (T-1)) \sum_{t=1}^{T-1}\frac{1}{t}\\
		&\leq 8 |\mS| \log\pi + 8 |\mS| \log(\pi (T-1)) (1 + \log (T-1))\\
		&\leq 8 |\mS| \log(\pi T) (2 + \log T).
	\end{align*}
	Here, the first inequality arises from Lemma~\ref{lemma:shrinking_difference}. The second and third inequalities are from bounding $t$ and from observing that $\sum_{t=1}^{T-1}(1/t) \leq 1 + \log(T-1)$. 
	
	\noindent Now for the function $h'(T)$, we can calculate
	\begin{align*}
		\sum_{t=2}^{T}\|p_t - q_t\| &\leq \sum_{t=2}^{T}4\frac{\sqrt{|\mS|\log(\pi ({t-1}))}}{\sqrt{{t-1}}}\\
		&\leq 4 |\mS|^{\frac{1}{2}} \sqrt{\log(\pi T)} \sum_{t=2}^{T}\frac{1}{\sqrt{{t-1}}}\\
		&\leq 8 |\mS|^{\frac{1}{2}} \sqrt{\log(\pi T)} \sqrt{T}.
	\end{align*}
	Here, the first inequality arises from Lemma~\ref{lemma:shrinking_difference}. The second and third inequalities are from bounding $t$ and from observing that $\sum_{t=2}^{T}(1/\sqrt{{t-1}}) \leq 2\sqrt{T-1} \leq 2 \sqrt{T}$. 
	
	\noindent \textbf{Kernel based ambiguity sets}.
	Calculating the function $h(T)$, we have
	\begin{align*}
		\frac{1}{2} \sum_{t=1}^{T}\|p_t - q_t\|_2^2 &\leq {\frac{1}{2}\left(2+ \frac{4\sqrt{C}}{\lambda}\right)^2}+  \sum_{t=2}^{T}\frac{32 C}{\lambda^2 ({t-1})}\log(\pi ({t-1})/\sqrt{6\delta})\\
		&\hspace{-10mm}\leq {\frac{1}{2}\left(2+ \frac{4\sqrt{C}}{\lambda}\right)^2}+ \frac{32C}{\lambda^2 } \log(\pi T/\sqrt{6\delta})\sum_{t=2}^{T}\frac{1}{t-1}\\
		&\hspace{-10mm}\leq {\frac{1}{2}\left(2+ \frac{4\sqrt{C}}{\lambda}\right)^2}+ \frac{32C}{\lambda^2 } \log(\pi T/\sqrt{6\delta})(1 + \log T).
	\end{align*}
	Here, the first inequality arises from Lemma~\ref{lemma:kernel_shrinking_difference}. The second and third inequalities are from bounding $t$ and from observing that $\sum_{t=1}^{T-1}(1/t) \leq 1 + \log(T-1) \leq 1 + \log T$. 
	
	Now, for the function $h'(T)$, we have 
	\begin{align*}
		\sum_{t=2}^{T}\|p_t - q_t\|_2 &\leq  \sum_{t=2}^{T}\frac{8}{\lambda}\frac{\sqrt{C}}{\sqrt{t-1}}\sqrt{ \log(\pi ({t-1})/\sqrt{6\delta})}\\
		&\leq \frac{8\sqrt{C}}{\lambda} \sqrt{ \log(\pi T/\sqrt{6\delta})} \sum_{t=2}^{T}\frac{1}{\sqrt{t-1}}\\
		&\leq \frac{16\sqrt{C}}{\lambda} \sqrt{T \log(\pi T/\sqrt{6\delta})}. 	
	\end{align*}
	Here, the first inequality arises from Lemma~\ref{lemma:kernel_shrinking_difference}. The second and third inequalities are from bounding $t$ and from observing that $\sum_{t=1}^{T-1}(1/\sqrt{t-1}) \leq 2\sqrt{T-1} \leq 2 \sqrt{T}$. 
	
	\noindent \textbf{$\ell_2$-norm ambiguity sets}.
	Calculating the function $h(T)$, we have
	\begin{align*}
		\frac{1}{2} \sum_{t=1}^{T}\|p_t - q_t\|^2 &\leq \frac{1}{2}\left(4\sqrt{|\mS| \log (\pi / \sqrt{3 \delta})}\right)^2 + \frac{1}{2} \sum_{t=2}^{T}16 |\mS| \frac{\log(\pi ({t-1})/\sqrt{3 \delta})}{t-1}\\
		&\leq 8 |\mS| \log (\pi / \sqrt{3 \delta}) + 8 |\mS| \log(\pi T/\sqrt{3\delta}) \sum_{t=2}^{T}\frac{1}{t-1}\\
		&\leq 8 |\mS| \log(\pi T/\sqrt{3\delta}) (2 + \log T).
	\end{align*}
	Here, the first inequality arises from Lemma~\ref{lemma:l2_shrinking_difference}. The second and third inequalities are proven similar to the case of $h(T)$ for the interval sets. 
	
	\noindent Now for the function $h'(T)$, we can calculate
	\begin{align*}
		\sum_{t=2}^{T}\|p_t - q_t\| &\leq  \sum_{t=2}^{T}4\sqrt{\frac{|\mS| \log (\pi {(t-1)} / \sqrt{3 \delta})}{t-1}}\\
		&\leq 4\sqrt{|\mS| \log (\pi T / \sqrt{3 \delta})} \sum_{t=2}^{T}\sqrt{\frac{1}{t-1}}\\
		&\leq 8\sqrt{|\mS| T\log (\pi T / \sqrt{3 \delta})}. 
	\end{align*}
	Here, the first inequality arises from Lemma~\ref{lemma:l2_shrinking_difference}. The second and third inequalities are proven similar to the case of $h'(T)$ for the interval sets.  \Halmos


\subsection{Numerical Experiments}

In this section, we provide the details
for the numerical experiments
conducted in Section~\ref{sec:num_res}.
We also provide an additional set of experiments
on an optimal routing problem to further illustrate our algorithms. 

\subsubsection{Benchmark Instances}

All mixed-integer linear optimization problems used in our numerical experiments can be found in Table \ref{Table:overview_mip}. The entries show the number of variables and constraints for each problem. 
The same holds for the quadratic problems listed in Table \ref{Table:overview_miqp}.

\begin{table}[hbt]
	\centering
	\begin{tabular}{lrrrrr}
		\toprule
		Name &  \multicolumn{4}{c}{Variables} & Constraints \\
		\midrule
		& All & Bin. & Int. & Cont. & \\ \midrule
		blend2 &  353 & 239 & 25 & 89 & 274 \\
		flugpl & 18 & 0 & 11 & 7 & 19 \\
		gr4x6 &  48 & 24 & 0 & 24 & 34 \\
		neos-1430701 	&312& 	156& 	0 &	156 &	668 \\
		noswot& 	128 &	75 	&25 &	28 &	182\\
		prod1& 	250 &	149 &	0 	&101& 	208 \\
		prod2&	301 &	200 &	0 &	101& 	211 \\ 
		ran13x13 &		338 &	169 &	0 	&169 &	195 \\
		supportcase14& 	304 &	304 &	0& 	0 &	234 \\
		supportcase16& 	319 &	319 &	0& 	0 &	130 \\
		beavma& 390 &	195 &	0 &	195 &	372\\
		k16x240b& 480& 	240& 	0 &240 &	256 \\
		neos-3610040 &430 &	85 &	0 &	345 &	335  \\
		neos-3611689 &421 	&88 &	0 	&333& 	323 \\
		timtab1CUTS & 	397& 	77 	&94 &	226 &	371 \\
		\bottomrule
	\end{tabular}
	\caption{\label{Table:overview_mip}Overview of MIP instances}
\end{table}

\begin{table}[htb]
	\centering
	\begin{tabular}{lrrrrr}
				\toprule
		Name &  \multicolumn{3}{c}{Variables} & \multicolumn{2}{c}{Constraints} \\
		\midrule
		& All & Bin.  & Cont. & All &  Quadr. \\ \midrule
		7579 & 300 & 100 & 200 & 203 & 1 \\
		10001 & 485 & 426  & 59 & 296 & 1 \\
		10002 & 485 & 426  & 59 & 296 & 1 \\
		10004 & 1058 & 999  & 59 & 867 & 1 \\
		10010 & 269 & 262  & 7 & 147 & 1 \\
		10003 & 1058 & 999  & 59 & 867 & 1\\
		10008 & 845 & 713  & 132 & 416 & 1\\
		10009 & 605 & 473  & 132 & 246 & 1\\
		10011 & 1390 & 1258 & 132 & 873 & 1 \\ 
		10012 & 967 & 835 & 132 & 538 & 1 \\
		\bottomrule
	\end{tabular}
	\caption{Overview of MIQP instances}\label{Table:overview_miqp}
\end{table}